\documentclass[reqno,oneside,a4paper,12pt]{amsart}

\usepackage{pstricks}
\usepackage{verbatim}
\usepackage{amssymb}
\usepackage{mathrsfs}
\usepackage{dsfont}
\usepackage{a4wide}
\usepackage{graphicx}

\AtBeginDocument{%
}


\newtheorem{lem}{Lemma}[section]
\newtheorem{thm}[lem]{Theorem}
\newtheorem{cor}[lem]{Corollary}
\newtheorem{prop}[lem]{Proposition}

\theoremstyle{definition}
\newtheorem{defn}[lem]{Definition}
\newtheorem*{notation}{Notation}
\newtheorem{rem}[lem]{Remark}

\newcommand{\stab}{\mathrm{stab}}
\newcommand{\bu}{\mathbb{U}}
\newcommand{\conv}{\mathrm{conv}}
\newcommand{\taua}{\tau_{\alpha\vphantom{/2}}}
\newcommand{\tauah}{\tau_{\alpha/2}^{-1/2}}
\newcommand{\Hom}{\mathrm{Hom}}

\renewcommand{\a}{\alpha}

\renewcommand{\S}{\Sigma}
\renewcommand{\t}{\theta}
\newcommand{\la}{\lambda}

\newcommand{\at}{\tilde{\alpha}^{\vee}}
\renewcommand{\o}{\omega}
\renewcommand{\O}{\Omega}
\newcommand{\s}{\sigma}
\newcommand{\ca}{\mathcal{A}}
\newcommand{\cb}{\mathcal{B}}
\newcommand{\cc}{\mathcal{C}}

\newcommand{\ch}{\mathcal{H}}

\newcommand{\cs}{\mathcal{S}}

\newcommand{\bc}{\mathbb{C}}
\newcommand{\bn}{\mathbb{N}}
\newcommand{\br}{\mathbb{R}}
\newcommand{\bt}{\mathbb{T}}
\newcommand{\bz}{\mathbb{Z}}

\newcommand{\G}{\Gamma}

\newcommand{\sca}{\mathscr{A}}

\newcommand{\lan}{\langle}
\newcommand{\ran}{\rangle}

\newcommand{\Aut}{\mathrm{Aut}}

\newcommand{\scx}{\mathscr{X}}

\newcommand{\ts}{\textsection}

\newcommand{\cha}{\alpha^{\vee}}
\newcommand{\chR}{R^{\vee}}

\DeclareMathAlphabet{\mathpzc}{OT1}{pzc}{m}{it}

\newcommand{\e}{\epsilon}

\numberwithin{equation}{section} \numberwithin{figure}{section}

\begin{document}

\title{Isotropic Random Walks on Affine Buildings}
\author{James Parkinson}

\date{\today}

\begin{abstract} Recently, Cartwright and Woess \cite{cw} provided a detailed
analysis of \textit{isotropic} random walks on the vertices of
thick affine buildings of type $\tilde{A}_n$. Their results
generalise results of Sawyer \cite{ss} where homogeneous trees
are studied (these are $\tilde{A}_1$ buildings), and Lindlbauer
and Voit \cite{hyper}, where $\tilde{A}_2$ buildings are studied.
In this paper we apply techniques of spherical harmonic analysis
to prove a local limit theorem, a rate of escape theorem, and a
central limit theorem for isotropic random walks on arbitrary
thick regular affine buildings of irreducible type, thus providing
a broad generalisation of the $\tilde{A}_n$ case.
\end{abstract}

\maketitle

\section*{Introduction}

Let $\scx$ be a thick locally finite regular affine building of
irreducible type. By \textit{regular} we mean that the number of
chambers containing a panel depends only on the cotype of the
panel, and by \textit{thick} we mean that this number is always at
least $3$. The simplest example of such a building is a
homogeneous tree with degree $q+1\geq3$, where the chambers are
the edges of the graph. In this case, Sawyer \cite{ss} studied isotropic
random walks $(Z_k)_{k\geq0}$ on the vertices of $\scx$, meaning that the transition probabilities
$p(x,y)=\mathbb{P}(Z_{k+1}=y\mid Z_k=x)$ depend only on the graph
distance $d(x,y)$ between $x$ and ~$y$. To motivate our results,
let us briefly describe these random walks on trees.

Let $V$ be the vertex set of the tree, and for each $x\in V$ and
$k\in\bn=\{0,1,\ldots\}$, write $V_k(x)$ for the set of all $y\in
V$ such that $d(x,y)=k$. It is easily seen that the cardinalities
$|V_k(x)|$ are independent of the particular $x\in V$, and we
write $N_k$ for this value. For each $k\in\bn$ there is a natural
operator $A_k$ acting on the space of all functions $f:V\to\bc$,
where for each $x\in V$, $(A_kf)(x)$ is the average value of $f$
over $V_{k}(x)$. The operator $A_k$ may be regarded as the
transition operator of the isotropic random walk with matrix
$(p_k(x,y))_{x,y\in V}$, where $p_k(x,y)=\frac{1}{N_k}$ if $y\in
V_k(x)$ and $p_k(x,y)=0$ otherwise. Indeed it is easily seen that
a random walk on $V$ is isotropic if and only if it has a
transition operator of the form
\begin{align*}
A=\sum_{k\in\bn}a_kA_k
\end{align*}
where $a_k\geq0$ for all $k\in\bn$ and $\sum_{k\in\bn}a_k=1$.

The linear span over $\bc$ of the operators $\{A_k\}_{k\in\bn}$
is a commutative algebra $\sca$ with a rich theory of harmonic
analysis (see \cite{FN}). In particular, the algebra
homomorphisms $h:\sca\to\bc$ may be explicitly described, and
local limit theorems, central limit theorems, and rate of escape
theorems can be proved as applications.

Now consider a regular affine building $\scx$ of irreducible
type. Thus $\scx$ may be regarded as a simplicial complex made by
`gluing together' many copies of a given \textit{Coxeter
complex}, each Coxeter complex called an \textit{apartment}
of~$\scx$ (these are regular tessellations of Euclidean space by
simplices). There is an irreducible (but not necessarily reduced)
root system $R$ associated to $\scx$, as described in
Section~\ref{section:building}, and the \textit{coweight
lattice}~$P$ of~$R$ is a subset of the vertex set of the standard
Coxeter complex~$\S$. We consider random walks on a related
subset~$V_P$ of the vertices of~$\scx$, which in most cases is
the set of all special vertices of~$\scx$.

Let~$P^+$ be the set of dominant coweights of~$R$ (relative to
some fixed base). For each $x\in V_P$ there is a natural
partition of $V_P$ into subsets $V_{\la}(x)$, $\la\in P^+$, as
described in Definition~\ref{defn:V}. Roughly speaking, $y\in
V_{\la}(x)$ means that there exists an apartment $\ca$ containing
$x$ and~$y$ and a `suitable' isomorphism $\psi:\ca\to\S$ such
that $\psi(x)=0$ and $\psi(y)=\la\in P^+$ (in other words, $y$ is
\textit{in position $\la$ from $x$}). It is shown in
\cite[Theorem~5.15]{p} that for all $\la\in P^+$ the cardinality
of the set $V_{\la}(x)$ is independent of the particular $x\in
V_P$, and we write $N_{\la}$ for this value. Following the tree
case, for each $\la\in P^+$ let $A_{\la}$ be the operator acting
on the space of functions $f:V_P\to\bc$ with $(A_{\la}f)(x)$
being the average value of $f$ over $V_{\la}(x)$. The linear span
of these operators over $\bc$ is a commutative algebra~$\sca$,
which has been studied extensively in~\cite{p}. As shown in \cite{p2}, the algebra
homomorphisms $h:\sca\to\bc$ may be explicitly described both in
terms of the \textit{Macdonald spherical functions} and in terms
of an integral over the \textit{boundary} of~$\scx$.

We call a random walk $(Z_k)_{k\in\bn}$ on $V_P$ with transition
matrix $(p(x,y))_{x,y\in V_P}$ \textit{isotropic} if
$p(x,y)=p(x',y')$ whenever $y\in V_{\la}(x)$ and $y'\in
V_{\la}(x')$ for some $\la\in P^+$. As in the tree case, the
operators $A_{\la}$ may be regarded as the transition operators
of isotropic random walks with matrices $(p_{\la}(x,y))_{x,y\in
V_P}$, where $p_{\la}(x,y)=\frac{1}{N_{\la}}$ if $y\in
V_{\la}(x)$ and $p_{\la}(x,y)=0$ otherwise. It is easily seen
that a random walk on $V_P$ is isotropic if and only if it has a
transition operator of the form
\begin{align}\label{Adefn}
A=\sum_{{\la}\in P^+}a_{\la}A_{\la}
\end{align}
where $a_{\la}\geq0$ for all ${\la}\in P^+$ and $\sum_{\la\in P^+
}a_{\la}=1$.

In this paper we apply the spherical harmonic analysis associated
to the algebra $\sca$ to prove a local limit theorem, a central
limit theorem, and a rate of escape theorem for isotropic random
walks on~$V_P$. These results generalise the results in \cite{cw}
where $\tilde{A}_n$ buildings are studied (which in turn
generalise the corresponding results for homogeneous trees). Our
results may also be viewed as `building analogues' of well known
results concerning random walks on semisimple Lie groups (see
\cite{bougerol} for example).

Let us briefly outline the structure of this paper. In
Section~\ref{section:1} we give a summary of some background
material, mostly from \cite{p} and~\cite{p2}. This section
includes a discussion of root systems, Coxeter complexes and
buildings, the algebra $\sca$, and the spherical harmonic
analysis associated to this algebra. The main sections of this
paper are Sections~\ref{section:2}, \ref{section:3}
and~\ref{section:4}. In Section~\ref{section:2} we give our local
limit theorem for isotropic random walks on~$V_P$, describing the
asymptotic behaviour of the $k$-step transition probabilities
$p^{(k)}(x,y)=\mathbb{P}(Z_k=y\mid Z_0=x)$. We also give
necessary and sufficient conditions for irreducibility and
aperiodicity of the random walk, and in Remark~\ref{remark:2} we
outline some applications of our local limit theorem to random
walks on groups acting on buildings. In Section~\ref{section:3}
we prove our rate of escape theorem. For each $k\in\bn$, let
$\nu_k\in P^+$ be such that $Z_k\in V_{\nu_k}(Z_0)$. We show
that, with probability~$1$, the vector $\frac{1}{k}\nu_k$
converges to a vector $\gamma$ in the underlying vector space of
the root system~$R$. We apply our local limit theorem to show
that each component of $\gamma$ (relative to a set of fundamental
coweights of~$P$) is strictly positive. In
Section~\ref{section:4} we prove our central limit theorem,
showing that there is a positive definite matrix~$\G$ such that,
with $\gamma$ as above, the vector
$\frac{1}{\sqrt{k}}(\nu_k-k\gamma)$ tends in distribution to the
normal distribution $N(0,\G)$. In Appendix~\ref{appendix:1} we
determine the algebra homomorphisms $h:\sca\to\bc$ which are
bounded (generalising \cite[Theorem~4.7.1]{macsph}).

\section{Affine Buildings and Spherical Harmonic Analysis}\label{section:1}

\subsection{Root Systems and Weyl Groups}\label{rtsys} Root systems play a significant role in this
work. We fix the following notations and conventions, generally
following \cite[Chapter~VI]{bourbaki}.

Let $R$ be an irreducible, but not necessarily reduced, root
system in a real vector space~$E$ with inner product
$\lan\cdot,\cdot\ran$. The \textit{rank} of $R$ is $n$, the
dimension of $E$. Let $I_0=\{1,2,\ldots,n\}$ and
$I=\{0,1,2,\ldots,n\}$, and let $B=\{\alpha_i\mid i\in I_0\}$ be a
fixed base of $R$. Write $R^+$ for the set of positive roots
(relative to $B$), and let $R^{\vee}=\{\alpha^{\vee}\mid\alpha\in
R\}$ be the \textit{dual} root system of~$R$, where for each
$\alpha\in R$ we write
$\alpha^{\vee}=\frac{2\alpha}{\lan\alpha,\alpha\ran}$. Since $R$
is irreducible, by \cite[VI, \ts1, No.8,
Proposition~25]{bourbaki} there is a unique \textit{highest root}
\begin{align}\label{eq:m}
\tilde{\alpha}=\sum_{i\in I_0}m_i\alpha_i
\end{align}
with the property that if $\beta=\sum_{i\in I_0}k_i\alpha_i\in R$
then $m_i\geq k_i$ for each $i\in I_0$.

For each $i\in I_0$ define $\la_i\in E$ by $\lan
\la_i,\alpha_j\ran=\delta_{i,j}$ for all $j\in I_0$. The elements
$\{\la_i\}_{i\in I_0}$ are called the \textit{fundamental
coweights of $R$}. The \textit{coweight lattice of $R$} is
$P=\sum_{i\in I_0}\bz\la_i$, and elements $\la\in P$ are called
\textit{coweights of $R$}. A coweight $\la\in P$ is said to be
\textit{dominant} if $\lan\la,\alpha_i\ran\geq0$ for all $i\in
I_0$, and we write $P^+$ for the set of all dominant coweights.
Let $Q=\sum_{\alpha\in R}\bz\alpha^{\vee}$ be the \textit{coroot
lattice of $R$}, and let $Q^+=\sum_{\alpha\in
R^+}\bn\alpha^{\vee}$. Note that $Q\subseteq P$, and by \cite[VI, \ts1, No.8,
Proposition~25]{bourbaki} we have $\tilde{\alpha}^{\vee}\in P^+$.

For each $\alpha\in R$ and $k\in\bz$ let $H_{\alpha;k}=\{x\in
E\mid\lan x,\alpha\ran=k\}$. We call these sets \textit{affine
hyperplanes}, or simply \textit{hyperplanes}. For each $\alpha\in
R$ and $k\in \bz$ let $s_{\alpha;k}$ denote the orthogonal
reflection in $H_{\alpha;k}$. Thus $s_{\alpha;k}(x)=x-(\lan
x,\alpha\ran-k)\alpha^{\vee}$ for all $x\in E$. Write
$s_{\alpha}$ in place of $s_{\alpha;0}$, $s_i$ in place of
$s_{\alpha_i}$ (for $i\in I_0$), and let
$s_0=s_{\tilde{\alpha};1}$. Let $W_0=W_0(R)$ be the \textit{Weyl
group of $R$}, and let $W=W(R)$ be the \textit{affine Weyl group
of $R$}. Thus $W_0$ is the subgroup of $\mathrm{GL}(E)$ generated
by $S_0=\{s_i\}_{i\in I_0}$, and $W$ is the subgroup of
$\mathrm{Aff}(E)$ generated by $S=\{s_i\}_{i\in I}$. Both
$(W_0,S_0)$ and $(W,S)$ are Coxeter systems, and clearly
$W_0\leq W$. Given $w\in W$, we define the \textit{length}
$\ell(w)$ of $w$ to be smallest $k\in\bn$ such that
$w=s_{i_1}\ldots s_{i_k}$, with $i_1,\ldots,i_k\in I$.

The \textit{extended affine Weyl group of $R$} is
$\tilde{W}(R)=\tilde{W}=W_0\ltimes P$. Since $W\cong W_0\ltimes
Q$ (see \cite[VI, \ts2, No.1, Proposition~1]{bourbaki}) and
$Q\subseteq P$, we may regard $W$ as a subgroup of~$\tilde{W}$.
Note that $\tilde{W}$ contains all translations by elements
of~$P$, while $W$ only contains those translations by elements
of~$Q$.

\begin{rem}\label{rem:BCn}
We make the following comments for those readers not so familiar
with the non-reduced root systems. For each $n\geq1$ there is
exactly one irreducible non-reduced root system (up to
isomorphism) of rank $n$, denoted by $BC_n$. To describe this
root system we may take $E=\br^n$ with the usual inner product,
and let $R$ consist of the vectors $\pm e_i,\pm 2e_i$ and $\pm
e_j\pm e_k$ for $1\leq i\leq n$ and $1\leq j<k\leq n$. Let
$\alpha_i=e_i-e_{i+1}$ for $1\leq i<n$ and $\alpha_n=e_n$. Then
$B=\{\alpha_i\}_{i\in I_0}$ is a base of $R$, and $R^+$ consists
of the vectors $e_i,2e_i$ and $e_j\pm e_k$ for $1\leq i\leq n$
and $1\leq j<k\leq n$. The fundamental coweights are
$\la_i=e_1+\cdots+e_i$ for each $1\leq i\leq n$. Note that
$\chR=R$ and $Q=P$. The subsystem $R_1=\{\alpha\in R\mid
2\alpha\notin R\}$ is a root system of type $C_n$, and the
subsystem $R_2=\{\alpha\in R\mid \frac{1}{2}\alpha\notin R\}$ is
a root system of type $B_n$ (with the convention that
$C_1=B_1=A_1$). We have $Q(R)=P(R)=Q(R_1)\subset P(R_1)$ (with strict inclusion),
and so $W(R)=W(R_1)$ but $\tilde{W}(R)\neq\tilde{W}(R_1)$.
\end{rem}

\subsection{The Coxeter Complex}\label{type}
There is a natural \textit{geometric realisation}~$\S=\S(R)$ of
the \textit{Coxeter complex of $W=W(R)$}. Let $\ch$ denote the
family of the hyperplanes $H_{\alpha;k}$, $\alpha\in R$,
$k\in\bz$, and define \textit{chambers} of $\S$ to be open
connected components of $E\backslash\bigcup_{H\in \ch}H$. Since
$R$ is irreducible each chamber is an open (geometric) simplex
\cite[V, \ts3, No.9, Proposition~8]{bourbaki}. We call the
extreme points of the closure of chambers \textit{vertices} of
$\S$, and we write $V(\S)$ for the set of all vertices of $\S$.

The set $P$ of coweights of $R$ is a subset of $V(\S)$, and we
call elements of $P$ the \textit{good vertices of $\S$}. When $R$
is reduced, $P$ is the set of more familiar \textit{special
vertices} of $\S$ \cite[VI, \ts2, No.2, Proposition~3]{bourbaki}.

The choice of the base $B=\{\alpha_i\}_{i=1}^n$ gives a natural
choice of a \textit{fundamental chamber}
\begin{align}
\label{fundamentalalcove}C_0=\{x\in E\mid\lan
x,\alpha_i\ran>0\textrm{ for all $i\in I_0$ and }\lan
x,\tilde{\alpha}\ran<1\}\,.
\end{align}
In the notation of~(\ref{eq:m}), the vertices of $C_0$ are the
points $\{0\}\cup\{\la_i/m_i\}_{i\in I_0}$ \cite[VI, \ts2,
No.2]{bourbaki}. There is a natural simplicial complex structure
on~$\S$ with maximal simplices being the vertex sets of chambers
of~$\S$, and simplices being subsets of the maximal simplices. We
define $\tau:V(\S)\to I$ to be the unique labelling of $\S$ (as a
simplicial complex) such that $\tau(0)=0$ and $\tau(\la_i/m_i)=i$
for each $i\in I_0$.

We write $I_{P}=\{\tau(\la)\mid\la\in P\}\subseteq I$. Let
$\{m_i\}_{i\in I_0}$ be as in (\ref{eq:m}), and define~$m_0=1$.
We have $I_P=\{i\in I\mid m_i=1\}$, which shows that $0\in I_P$
for all root systems, and that $I_P=\{0\}$ if $R$ is non-reduced
\cite[Lemma~4.3]{p}. This also shows that in the non-reduced
case, and only in the non-reduced case, there are special
vertices which are not good vertices.

We define the \textit{fundamental sector} of $\S$ to be the open
simplicial cone
\begin{align}
\label{fundamentalsector}\cs_0=\{x\in E\mid\lan
x,\alpha_i\ran>0\textrm{ for all $i\in I_0$}\}.
\end{align}
The \textit{sectors} of $\S$ are then the sets $\la+w\cs_0$,
where $w\in W_0$ and $\la\in P$ (equivalently, the sectors are the sets
$\tilde{w}\cs_0$, $\tilde{w}\in \tilde{W}$).

An \textit{automorphism of $\S$} is a
bijection $\psi$ of $E$ which maps chambers, and only chambers, to
chambers, with the property that chambers $C$ and $D$ are
adjacent if and only if $\psi(C)$ is adjacent to~$\psi(D)$. We
write $\Aut(\S)$ for the automorphism group of~$\S$. An automorphism $\psi$ of $\S$ is called \textit{type preserving}
if $\tau(v)=\tau(\psi(v))$ for all $v\in V(\S)$. By
\cite[Lemma~2.2]{ronan} $\psi\in\Aut(\S)$ is type preserving if
and only if~$\psi\in W$. Generally we have $W_0<W\leq
\tilde{W}\leq \Aut(\S)$ (with the possibility that $W<\tilde{W}$ and
$\tilde{W}<\Aut(\S)$).

\subsection{Buildings and Regularity}\label{section:building} Recall (\cite{brown}) that a \textit{building of type $W$} is a
nonempty simplicial complex $\scx$ which contains a family of
subcomplexes called \textit{apartments} such that
\begin{enumerate}
\item[$\mathrm{(i)}$] each apartment is isomorphic to the (simplicial) Coxeter
complex of $W$,
\item[$\mathrm{(ii)}$] given any two chambers of $\scx$ there is an
apartment containing both, and
\item[$\mathrm{(iii)}$] given any
two apartments $\ca$ and $\ca'$ that contain a common chamber,
there exists an isomorphism $\psi:\ca\to\ca'$ fixing $\ca\cap\ca'$
pointwise.
\end{enumerate}
Since $W$ is an affine Weyl group, $\scx$ is called an
\textit{affine building}.

It is an easy consequence of this definition that $\scx$ is a
labellable simplicial complex, and all the isomorphisms in the
above definition may be taken to be type preserving (this ensures
that the labellings of $\scx$ and $\S$ are \textit{compatible}).

Let $V$ and $\cc$ be the vertex and chamber sets of $\scx$,
respectively (with chambers being maximal simplices of $\scx$). Chambers $c$ and $d$ are declared to be
\textit{$i$-adjacent} (written $c\sim_i d$) if and only if either $c=d$, or if all the
vertices of $c$ and $d$ are the same except for those of
type~$i$.

Throughout this paper we assume that our buildings are
\begin{enumerate}
\item[(i)] \textit{locally finite}, meaning that $|I|<\infty$
and $|\{d\in\cc\mid d\sim_i c\}|<\infty$ for each $c\in\cc$ and
$i\in I$,
\item[(ii)] \textit{regular}, meaning that $|\{d\in\cc\mid d\sim_i c\}|$ is independent of $c\in\cc$ for
each $i\in I$, and
\item[(iii)] \textit{thick}, meaning that $|\{d\in\cc\mid d\sim_i
c\}|\geq3$ for each $c\in\cc$ and each $i\in I$.
\end{enumerate}

By \cite[Theorem~2.4]{p} we see that thickness and regularity are
intimately connected. Indeed, the only thick affine buildings of
irreducible type which may fail to be regular are those of dimension~1
(thus regularity is a very weak hypothesis).

Since $\scx$ is assumed to be regular, we may define numbers $q_i$,
$i\in I$, called the \textit{parameters} of the building, by
$q_i+1=|\{d\in\cc\mid d\sim_i c\}|$. These numbers satisfy
$q_j=q_i$ if $s_j=ws_iw^{-1}$ for some $w\in W$ (see
\cite[Corollary~2.2]{p}), and by thickness $q_i>1$ for all $i\in
I$. If $w=s_{i_1}\cdots s_{i_k}\in W$ is a reduced expression
(that is, $\ell(w)=k$) we define $q_w=q_{i_1}\cdots q_{i_k}$,
which is independent of the particular reduced expression for $w$
(see \cite[Proposition~2.1(i)]{p}).

To each locally finite regular affine building of irreducible
type we associate an irreducible root system~$R$ (depending on
the parameter system of the building) as follows (see
\cite[Appendix]{p}):
\begin{enumerate}
\item[(i)] If $\scx$ is a regular $\tilde{A}_1$ building with
$q_0=q_1$, then we take $R=A_1$ (these buildings are homogeneous
trees).
\item[(ii)] If $\scx$ is a regular $\tilde{A}_1$ building with $q_0\neq q_1$, then we
take $R=BC_1$ (these buildings are \textit{semi-homogeneous
trees}).
\item[(iii)] If $\scx$ is a regular $\tilde{C}_n$ building with
$n\geq2$ and $q_0=q_n$, then we take $R=C_n$.
\item[(iv)] If $\scx$ is a regular $\tilde{C}_n$ building with $n\geq2$ and $q_0\neq q_n$, then we take
$R=BC_n$.
\item[(v)] If $\scx$ is a regular building of type $\tilde{X}_n$, where
$X=A$ and $n\geq2$, or $X=B$ and $n\geq3$, or $X=D$ and $n\geq4$,
or $X=E$ and $n=6,7$ or $8$, or $X=F$ and $n=4$, or $X=G$ and
$n=2$, then we take $R=X_n$.
\end{enumerate}

The choices above are made to ensure that the coweight lattice
$P$ of $R$ preserves the parameter system of $\scx$ in the sense
that if $v\in V(\S)$ then $q_{\tau(v)}=q_{\tau(v+\la)}$ for all
$\la\in P$. Thus, for example,~(iv) above is motivated by the
general parameter system of a $\tilde{C}_n$ building (embodied in
the Coxeter graph):\newline
\begin{figure}[ht]
 \begin{center}
 \psset{xunit= 0.9 cm,yunit= 0.9 cm}
 \psset{origin={0,0}}
\vspace{0.3cm}

\rput(-3,0.3){$q_0$}\pscircle*(-3,0){2pt}
\rput(-2,0.3){$q_1$}\pscircle*(-2,0){2pt} \psline(-3,0)(-2,0)
\rput(-2.5,0.3){$4$} \rput(-1,0.3){$q_1$}\pscircle*(-1,0){2pt}
\psline(-2,0)(-1,0) \pscircle*(0.4,0){1.5pt}
\pscircle*(-0.4,0){1.5pt} \pscircle*(0,0){1.5pt}
\rput(1,0.3){$q_1$}\pscircle*(1,0){2pt}
\rput(2,0.3){$q_1$}\pscircle*(2,0){2pt} \psline(1,0)(3,0)
\rput(3,0.3){$q_n$}\pscircle*(3,0){2pt} \rput(2.5,0.3){$4$}
\end{center}
\caption{}\label{fig:fig1}
\end{figure}

\noindent(see \cite[Appendix]{p}). If we take $R=C_n$, then by
the definition of the type map (see Section~\ref{type}) and the
fact that $m_n=1$ (see (\ref{eq:m}) and
\cite[Plate~III]{bourbaki}) we have $\tau(\la_n)=n$. Thus in
general we have $q_{\tau(\la_n)}=q_n\neq q_0=q_{\tau(0)}$. If we
instead choose $R=BC_n$, then $P=Q$, and so $\tau(v+\la)=\tau(v)$
for all $v\in V(\S)$ and $\la\in P$ (and hence
$q_{\tau(v)}=q_{\tau(v+\la)}$).

\begin{defn}\label{goodbuilding} Let $\scx$ be a regular affine building with associated root system~$R$ and vertex
set~$V$. A vertex $x\in V$ is said to be \textit{good} if and only
if $\tau(x)\in I_P$ (recall that $I_P=\{\tau(\la)\mid\la\in P\}$).
Write $V_P$ for the set of all good vertices of~$\scx$.
\end{defn}

It is clear that $V_P$ is a subset of the more familiar
\textit{special vertices} of $\scx$. In fact if $R$ is reduced
then $V_P$ equals the set of special vertices. If $R$ is
non-reduced (so $R$ is of type $BC_n$ for some $n\geq1$), then
$V_P$ is the set of all type~$0$ vertices of~$\scx$ (whereas the
special vertices are those with types~$0$ or~$n$).

\subsection{The Algebra $\sca$}\label{uu} In this section we describe a commutative algebra
$\sca$ of \textit{vertex set averaging operators}. This algebra has been studied
in detail in~\cite{p}, where it is shown that $\sca$ is isomorphic to the center
of an appropriate affine Hecke algebra.

\begin{defn} Let $\ca$ be an apartment of $\scx$. An isomorphism
$\psi:\ca\to\S$ is called \textit{type rotating} if and only if
it is of the form $\psi=w\circ\psi_0$, where $\psi_0:\ca\to\S$ is
a type preserving isomorphism, and $w\in \tilde{W}$.
\end{defn}

\begin{defn}\label{defn:V} Given $x\in V_P$ and $\la\in P^+$, let
$V_{\la}(x)$ be the set of all $y\in V_P$ such that there exists
an apartment $\ca$ containing $x$ and $y$, and a type rotating
isomorphism $\psi:\ca\to\S$ such that $\psi(x)=0$ and
$\psi(y)=\la$. Equivalently, $y\in V_{\la}(x)$ if and only if
there exists an apartment $\ca$ containing $x$ and~$y$ and a type
rotating isomorphism $\psi:\ca\to\S$ such that $\psi(x)=0$ and
$\psi(y)\in W_0\la$.
\end{defn}

The requirement that $\psi$ is type rotating in
Definition~\ref{defn:V} ensures that $y\in V_{\la}(x)\cap
V_{\mu}(x)$ implies that $\la=\mu$. Indeed, in
\cite[Proposition~5.6]{p} we showed that for each $x\in V_P$,
$\{V_{\la}(x)\}_{\la\in P^+}$ forms a partition of~$V_P$.

\begin{rem}\label{rem:A1} To get a feel for Definition~\ref{defn:V} in a
special case, suppose that $\scx$ is a homogeneous tree with
degree $q+1$. Thus $R$ has type~$A_1$, and we may take
$R=\{\alpha,-\alpha\}$ where $\alpha=e_1-e_2$ (the underlying
vector space here is $E=\{\xi\in\br^2\mid\xi_1+\xi_2=0\}$).
Taking $B=\{\alpha\}$ we have $\la_1=\frac{\alpha}{2}$ and
$P^+=\{k\la_1\mid k\in\bn\}$. We have $V_P=V$ (all vertices are
`good'), and writing $V_k(x)$ in place of $V_{k\la_1}(x)$, we have
$$V_k(x)=\{y\in V\mid d(x,y)=k\},$$
where $d:V\times V\to\bn$ is the usual graph metric.

Note that in this example all isomorphisms $\psi:\ca\to\S$ where
$\ca$ is an apartment of $\scx$ are type rotating. To understand
why the type rotating hypothesis becomes important, suppose that
$\scx$ is a regular $\tilde{A}_2$ building, and take vertices
$x,y\in V_P$ with $y\in V_{\la_1}(x)$. Thus there exists an
apartment $\ca$ containing $x$ and $y$ and a type rotating
isomorphism $\psi:\ca\to\S$ with $\psi(x)=0$ and $\psi(y)=\la_1$.
The map $\varphi:\S\to\S$ given by $a_1\la_1+a_2\la_2\mapsto
a_2\la_1+a_1\la_2$ is an automorphism of $\S$, and so
$\varphi\circ\psi:\ca\to\S$ is an isomorphism (however it is
\textit{not} type rotating). Notice that
$(\varphi\circ\psi)(x)=0$ and $(\varphi\circ\psi)(y)=\la_2$, and
so if we drop the type rotating hypothesis in
Definition~\ref{defn:V} we would conclude that $y\in
V_{\la_1}(x)\cap V_{\la_2}(x)$.
\end{rem}

For $\la\in P$ let $\la^*=-w_0\la$, where $w_0$ is the unique
longest element of $W_0$. In \cite[Proposition~5.8]{p} we showed
that if $\la\in P^+$ then $\la^*\in P^+$, and that $y\in
V_{\la}(x)$ if and only if $x\in V_{\la^*}(y)$. Note that $*$ is
trivial unless $w_0\neq-1$, that is, unless $R=A_n, D_{2n+1}$ or
$E_6$ for some $n\geq2$ (see \cite[Plates~I-IX]{bourbaki}). For
example, the map $\varphi$ from Remark~\ref{rem:A1}
is~$\la\mapsto\la^*$.

In \cite[Theorem~5.15]{p} we showed that
$|V_{\la}(x)|=|V_{\la}(y)|$ for all $x,y\in V_P$ and $\la\in
P^+$, and we denote this common value by $N_{\la}$ (see~(\ref{eq:N}) for a
formula for~$N_{\la}$). For each $\la\in P^+$
define an operator $A_{\la}$, acting on the space of functions
$f:V_P\to\bc$, by
\begin{align*}
(A_{\la}f)(x)=\frac{1}{N_{\la}}\sum_{y\in
V_{\la}(x)}f(y)\quad\textrm{for all $x\in V_P$}
\end{align*}
(thus $(A_{\la}f)(x)$ is the average value of $f$ over the set
$V_{\la}(x)$). The linear span $\sca$ of $\{A_{\la}\}_{\la\in
P^+}$ over $\bc$ is a commutative algebra \cite[Theorem~5.24]{p}.

\begin{rem}\label{rem:A1(ii)} (i) In the situation of the first example of Remark~\ref{rem:A1}, writing
$N_{k}$ in place of $N_{k\la_1}$ we have $N_0=1$ and
$N_k=(q+1)q^{k-1}$ for $k\geq1$. In this case the operators
$A_{k}=A_{k\la_1}$ have been studied by many authors (see \cite[p.57]{FN}, \cite{ss} or \cite[\ts III.19.C]{woess}).
They satisfy the simple recurrence
$$A_kA_1=\frac{q}{q+1}A_{k+1}+\frac{1}{q+1}A_{k-1}\qquad\textrm{for
$k\geq1$},$$ although for general affine buildings such a formula
is not readily available.

(ii) Let $\sca_Q$ denote the linear span (over $\bc$) of
$\{A_{\la}\mid\la\in Q\cap P^+\}$. It is easily seen that~$\sca_Q$
is a subalgebra of~$\sca$. In the case when~$\scx$ is the
Bruhat-Tits building of a group~$G$ of $p$-adic type with maximal
compact subgroup $K$ (as in \cite[\ts2.4--2.7]{macsph}), $\sca_Q$
is isomorphic to $\mathscr{L}(G,K)$, the space of continuous
compactly supported bi-$K$-invariant functions on $G$.
\end{rem}

\subsection{Isotropic Random Walks}
As mentioned in the introduction, we call a random walk on $V_P$
with transition probability matrix $A=(p(x,y))_{x,y\in V_P}$
\textit{isotropic} if $p(x,y)=p(x',y')$ whenever $y\in
V_{\la}(x)$ and $y'\in V_{\la}(x')$ for some $\la\in P^+$. In
particular, each operator $A_{\la}$, $\la\in P^+$, represents an
isotropic random walk with transition matrix (also called
$A_{\la}$) given by $A_{\la}=(p_{\la}(x,y))_{x,y\in V_P}$, where
$p_{\la}(x,y)=N_{\la}^{-1}$ if $y\in V_{\la}(x)$ and
$p_{\la}(x,y)=0$ otherwise.

It is easily seen that a random walk is isotropic if and only if
its transition matrix (operator) $A$ is as in~(\ref{Adefn}). To
avoid triviality we always assume that $a_{\la}>0$ for at least
one $\la\neq0$ (so that $A$ is not the identity). In this paper
we will prove a local limit theorem, a rate of escape theorem,
and a central limit theorem for such random walks, generalising
the work of \cite{ss} (where homogeneous trees are studied) and
\cite{cw} (where $\tilde{A}_n$ buildings are studied). The main
techniques we use are those of \textit{spherical harmonic
analysis}, as recalled in the following sections. We note that
isotropic random walks on $\tilde{A}_2$ buildings have also been
studied by Lindlbauer and Voit \cite{hyper} where more hypergroup
oriented techniques are used (see \cite[\ts7]{p} for a discussion
of the hypergroups that arise in the setting of general affine
buildings).

In the case of Remark~\ref{rem:A1(ii)}(ii), the theorems we prove
in this paper can be translated into theorems concerning
probability measures on groups of $p$-adic type. We briefly
discuss this in Remark~\ref{remark:2}.

\subsection{The Algebra Homomorphisms $h:\sca\to\bc$}

Our proofs of the local limit theorem, rate of escape theorem and
central limit theorem rely heavily on two formulae for the
algebra homomorphisms $h:\sca\to\bc$. In this section we recall
these formulae from~\cite{p2}. The first formula is in terms of
the Macdonald spherical functions, and the second is in terms of
an integral over the \textit{boundary} of~$\scx$.

To simultaneously deal with the reduced and non-reduced cases we
introduce the following notation. Let $R_1=\{\alpha\in R\mid
2\alpha\notin R\}$, $R_2=\{\alpha\in R\mid
\frac{1}{2}{\alpha}\notin R\}$ and $R_3=R_1\cap R_2$. Notice that
$R_1=R_2=R_3=R$ if $R$ is reduced. For $\alpha\in R_2$, write
$q_{\alpha}=q_i$ if $|\alpha|=|\alpha_i|$ (if
$|\alpha|=|\alpha_i|$ then necessarily $\alpha\in R_2$). Since
$q_j=q_i$ whenever $s_j=ws_iw^{-1}$ for some $w\in W$ (see
\cite[Corollary~2.2]{p}), it follows that $q_i=q_j$ whenever
$|\alpha_i|=|\alpha_j|$, and so the definition of $q_{\alpha}$ is
unambiguous. Note that $R=R_3\cup(R_1\backslash R_3)\cup(R_2\backslash R_3)$
where the union is disjoint. Define a set of numbers
$\{\tau_{\alpha}\}_{\alpha\in R}$ related to the numbers
$\{q_{\alpha}\}_{\alpha\in R_2}$ by
\begin{align*}
\tau_{\alpha}=\begin{cases} q_{\alpha}&\textrm{if $\alpha\in
R_3$}\\ q_0&\textrm{if $\alpha\in R_1\backslash R_3$}\\
q_{\alpha}q_{0}^{-1}&\textrm{if $\alpha\in R_2\backslash
R_3$.}\end{cases}
\end{align*}
It is convenient to define $\tau_{\alpha}=1$ if $\alpha\notin R$.
Note that $\tau_{\alpha}=q_{\alpha}$ if $R$ is reduced (and many
subsequent formulae will simplify in this case).

If $u\in\Hom(P,\bc^{\times})$ we write $u^{\la}$ in place of
$u(\la)$. The homomorphism $r\in\Hom(P,\bc^{\times})$ defined by
\begin{align}\label{eq:r}
r^{\la}=\prod_{\alpha\in
R^+}\tau_{\alpha}^{\frac{1}{2}\lan\la,\alpha\ran}\qquad\textrm{for
all $\la\in P$}
\end{align}
plays an important role. By \cite[Proposition~1.5]{p2} and
\cite[Proposition~A.1]{p2} we have
\begin{align}\label{eq:N}
N_{\la}=N_{\la^*}=\frac{W_0(q^{-1})}{W_{0\la}(q^{-1})}r^{2\la}
\end{align}
where $W_{0\la}=\{w\in W_0\mid w\la=\la\}$ and where $X(q^{-1})=\sum_{w\in
X}q_w^{-1}$ for subsets $X\subset W_0$.

For $w\in W_0$ and $u\in\Hom(P,\bc^{\times})$ we write
$wu\in\Hom(P,\bc^{\times})$ for the homomorphism with
$(wu)^{\la}=u^{w\la}$ for all $\la\in P$. Following
\cite[Chapter~IV]{macsph}, for $\la\in P^+$ and
$u\in\Hom(P,\bc^{\times})$ we define the \textit{Macdonald
spherical function} $P_{\la}(u)$ by
\begin{align}\label{eq:mac}
P_{\la}(u)=\frac{r^{-\la}}{W_0(q^{-1})}\sum_{w\in
W_0}c(wu)u^{w\la}\quad\textrm{where}\quad c(u)=\prod_{\alpha\in
R^+}\frac{1-\taua^{-1}\tauah u^{-\cha}}{1-\tauah u^{-\cha}},
\end{align}
provided that the denominators of the $c(wu)$ functions do not
vanish. Since $P_{\la}(u)$ is a Laurent polynomial (see
\cite[(1.8)]{p2}), these \textit{singular} cases can be obtained
from the general formula by taking an appropriate limit (see
Lemma~\ref{later:refertothis} for one example).

For $u\in\Hom(P,\bc^{\times})$, let $h_u:\sca\to\bc$ be the
linear map with $h_u(A_{\la})=P_{\la}(u)$ for each $\la\in P^+$.
By \cite[Proposition~2.1]{p2} every algebra homomorphism
$h:\sca\to\bc$ is of the form $h=h_u$ for some $u\in
\Hom(P,\bc^{\times})$, and $h_{u'}=h_u$ if and only if $u'=wu$
for some $w\in W_0$. We call the formula
$h_u(A_{\la})=P_{\la}(u)$ the \textit{Macdonald formula for the
algebra homomorphisms $h:\sca\to\bc$}.

\begin{rem}\label{rem:A1(iii)} (i) In the situation of
homogeneous trees from Remarks~\ref{rem:A1}
and~\ref{rem:A1(ii)}(i), if $u\in\Hom(P,\bc^{\times})$, then
writing $z=u^{\la_1}\in\bc^{\times}$ we have
$$h_u(A_k)=\frac{q^{-k/2}}{1+q^{-1}}\left(\frac{1-q^{-1}z^{-2}}{1-z^{-2}}z^k+\frac{1-q^{-1}z^{2}}{1-z^2}z^{-k}\right)$$
provided that $z\neq\pm 1$ (with the values at $z=\pm1$ found by
taking appropriate limits). More generally, in the $\tilde{A}_n$ case the
functions $P_{\la}(u)$ are essentially the \textit{Hall-Littlewood polynomials} of
\cite{m2} (see \cite{C2}).

(ii) At times the $BC_n$ case (see Remark~\ref{rem:BCn}) requires
separate treatment. Recall the description of the parameter
system from Figure~\ref{fig:fig1}. For
$u\in\Hom(P,\bc^{\times})$, by writing $t_i=u^{e_i}$ for $1\leq
i\leq n$ (noting that in this case $e_i\in P$ for each $1\leq
i\leq n$), we have
$$c(u)=\bigg\{\prod_{i=1}^{n}\frac{(1-a^{-1}t_i^{-1})(1+b^{-1}t_i^{-1})}{1-t_i^{-2}}\bigg\}\bigg\{\prod_{1\leq j<k\leq
n}\frac{(1-q_1^{-1}t_j^{-1}t_k)(1-q_1^{-1}t_j^{-1}t_k^{-1})}{(1-t_j^{-1}t_k)(1-t_j^{-1}t_k^{-1})}\bigg\},$$
where $a=\sqrt{q_nq_0}$ and $b=\sqrt{q_n/q_0}$ (see
\cite[Section~5.2]{p2}).
\end{rem}

We now recall the second formula for the algebra homomorphisms
$h:\sca\to\bc$. A \textit{sector} of $\scx$ is a subcomplex
$\cs\subset\scx$ such that there exists an apartment $\ca$ with
$\cs\subset\ca$, and an isomorphism
$\psi:\ca\to\S$ such that $\psi(\cs)$ is a sector of~$\S$. The
\textit{base vertex} of $\cs$ is $\psi^{-1}(\la)$, where $\la\in
P$ is the base vertex of $\psi(\cs)$. If $\cs$ and $\cs'$ are
sectors of $\scx$ with $\cs'\subseteq\cs$, then we say that
$\cs'$ is a \textit{subsector} of $\cs$. The \textit{boundary}
$\O$ of $\scx$ is the set of equivalence classes of sectors,
where we declare two sectors to be equivalent if and only if they
contain a common subsector. Given $x\in V_P$ and $\o\in\O$ there
exists a unique sector, denoted $\cs^x(\o)$, in the class $\o$
with base vertex~$x$ \cite[Lemma~9.7]{ronan}. For each $x\in V_P$, $\o\in\O$ and $\la\in P^+$, the intersection
$V_{\la}(x)\cap\cs^x(\o)$ contains exactly one vertex, denoted
$v_{\la}^x(\o)$ (the reader is encouraged to draw a picture showing the vertices
$v_{\la}^x(\o)$ in the
$\tilde{A}_2$ case). By \cite[Theorem~3.4]{p2}, for each $\o\in\O$
and $x,y\in V_P$ there exists a coweight $h(x,y;\o)\in P$ such
that
\begin{align}\label{eq:horo}
v_{\mu}^x(\o)=v_{\mu-h(x,y;\o)}^y(\o)
\end{align}
for $\mu\in P^+$ with each $\lan\mu,\alpha_i\ran$, $i\in I_0$,
sufficiently large. Indeed, if $y\in V_{\la}(x)$ then
(\ref{eq:horo}) holds, for all $\o\in\O$, whenever
$\mu-\Pi_{\la}\subset P^+$ (see \cite[Theorem~3.6]{p2}). Here
$\Pi_{\la}\subset P$ is the \textit{saturated} set with highest
coweight~$\la$ relative to the partial order on $P$ given by
$\mu\preceq\la$ if and only if $\la-\mu\in Q^+$ (recall that
$Q^+$ is the $\bn$-span of $\{\alpha^{\vee}\mid\alpha\in R^+\}$).
We have
$$\Pi_{\la}=\{w\nu\mid \nu\in P^+,\nu\preceq\la,w\in W_0\}$$
(see \cite[Lemma~13.4B]{h2} for example). The vectors $h(x,y;\o)$
are generalisations of the so called \textit{horocycle numbers}
for homogeneous trees.

By \cite[Proposition~3.5]{p2}, for all $\o\in\O$ and all
$x,y,z\in V_P$ we have the \textit{cocycle relation}
\begin{align}\label{eq:cocycle}
h(x,y;\o)=h(x,z;\o)+h(z,y;\o).
\end{align}
Thus $h(x,x;\o)=0$ and
$h(x,y;\o)=-h(y,x;\o)$ for all $\o\in\O$ and all $x,y\in V_P$.

There is a natural topology on $\O$ (discussed in \cite{p2}) in
which for each $x\in V_P$ the sets $\O_x(y)=\{\o\in\O\mid
y\in\cs^y(\o)\}$, $y\in V_P$, form a basis of open and closed sets
(this topology is independent of the particular $x\in V_P$
chosen). For each $x\in V_P$ there is a unique regular Borel
probability measure $\nu_x$ on $\O$ such that
$\nu_x(\O_x(y))=N_{\la}^{-1}$ if $y\in V_{\la}(x)$. For $x,x'\in
V_P$ the measures $\nu_x$ and $\nu_{x'}$ are mutually absolutely
continuous with Radon-Nikodym derivative
$(d\nu_{x'}/d\nu_x)(\o)=r^{2h(x,x';\o)}$ (see
\cite[Theorem~3.17]{p2}).

The \textit{integral formula for the algebra homomorphisms
$h:\sca\to\bc$} is
\begin{align}\label{eq:twoforms}
P_{\la}(u)=h_u(A_{\la})=\int_{\O}(ur)^{h(x,y;\o)}d\nu_x(\o)
\end{align}~for any $x,y\in V_P$ with $y\in V_{\la}(x)$. Equality of the Macdonald and integral
formulae is non-trivial, and is proved in \cite[Corollary~3.23 and
Theorem~6.2]{p2}.

\subsection{The Plancherel measure}

The \textit{Plancherel measure} of $\sca$ is instrumental in our proof
of the local limit theorem. In this section we recall some details about the
Plancherel measure and the $\ell^2$-spectrum of $\sca$ from \cite{p2} (see also
\cite{macsph}).

It is easy to see that each $A\in\sca$ maps $\ell^2(V_P)$ into
itself, and for $\la\in P^+$ and $f\in\ell^2(V_P)$ we have
$\|A_{\la}f\|_2\leq\|f\|_2$ (see \cite[Lemma~4.1]{C2} for a proof
in a similar context). So we may regard $\sca$ as a subalgebra of
the $C^*$-algebra $\mathscr{L}(\ell^2(V_P))$ of bounded linear
operators on $\ell^2(V_P)$. The facts that $y\in V_{\la}(x)$ if
and only if $x\in V_{\la^*}(y)$, and $N_{\la^*}=N_{\la}$, imply
that $A_{\la}^*=A_{\la^*}$, and so the adjoint $A^*$ of any
$A\in\sca$ is also in~$\sca$.

Let $\sca_2$ denote the completion of $\sca$ with respect to
$\|\cdot\|$, the $\ell^2$-operator norm. So $\sca_2$ is a
commutative $C^*$-algebra. The algebra homomorphisms
$h:\sca_2\to\bc$ are precisely the extensions $h=\tilde{h}_u$ of
those algebra homomorphisms $h_u:\sca\to\bc$ which are continuous
with respect to the $\ell^2$-operator norm. Let us describe the
latter homomorphisms.

The analysis here splits into two cases. Following
\cite[Chapter~V]{macsph} we call the situation where
$\tau_{\alpha}\geq1$ for all $\alpha\in R$ the \textit{standard
case}, and the situation where $\tau_{\alpha}<1$ for some
$\alpha\in R$ the \textit{exceptional case} (the use of the word
``exceptional'' here is unrelated to the so called
\textit{exceptional root systems}). It is immediate from the
definition of the numbers $\tau_{\alpha}$ that the exceptional
case occurs exactly when $R=BC_n$ for some $n\geq1$ and $q_n<q_0$
(see \cite[Lemma~5.1]{p2}). In particular, if $R$ is reduced
then we are in the standard case.

Let us consider the standard case first. Let
$$\bu=\{u\in\Hom(P,\bc^{\times}): |u^{\la}|=1\textrm{ for all
$\la\in P$}\}.$$ In the standard case the algebra homomorphism
$h_u:\sca\to\bc$ is continuous with respect to the
$\ell^2$-operator norm if and only if $u\in\bu$ (see
\cite[Corollary~5.4]{p2}). If $h=\tilde{h}_u$, $u\in\bu$, we
write $\widehat{A}(u)=h(A)$ for $A\in\sca_2$. In particular,
$\widehat{A}_{\la}(u)=P_{\la}(u)$ for $u\in\bu$.

In the standard case, let $\pi$ be the measure on $\bu$ given by
$d\pi(u)=\frac{W_0(q^{-1})}{|W_0|}|c(u)|^{-2}du$, where $du$ is
normalised Haar measure on $\mathbb{U}$ (note that in \cite{p2}
we write $\pi_0$ instead of $\pi$). Then for $A\in\sca_2$ we have
\begin{align*}
(A\delta_y)(x)=\int_{\mathbb{U}}\widehat{A}(u)\overline{P_{\la}(u)}d\pi(u)\qquad\textrm{whenever
$y\in V_{\la}(x)$}
\end{align*}
where $\delta_y(x)=1$ if $x=y$ and $\delta_y(x)=0$ otherwise (see
\cite[Theorem~5.2 and Corollary~5.5]{p2}). The measure $\pi$ is
essentially the Plancherel measure of $\sca$ (more precisely, the
Plancherel measure is the image of the measure $\pi$ under the
homeomorphism $\varpi:\bu/W_0\to\Hom(\sca_2,\bc)$,
$u\mapsto\tilde{h}_u$).

Let us consider the exceptional case, and so $R=BC_n$ for some
$n\geq1$ and $q_n<q_0$. For $u\in \Hom(P,\bc^{\times})$, recall
the definition of the numbers $t_i=t_i(u)$, $1\leq i\leq n$, from
Remark~\ref{rem:A1(iii)}(ii). We use the isomorphism
$\bu\to\bt^n$, $u\mapsto(t_1,\ldots,t_n)$ to identify $\bu$ with
$\bt^n$ (here $\bt=\{t\in\bc:|t|=1\}$). Define $\bu'=\{-b\}\times
\bt^{n-1}$, and write $U=\bu\cup\bu'$ (recall from
Remark~\ref{rem:A1(iii)}(ii) that $b=\sqrt{q_n/q_0}$). Let $dt=dt_1\cdots dt_n$, where $dt_i$ is normalised Haar measure
on $\bt$. Let $\phi_0(u)=c(u)c(u^{-1})$, and let
$$\phi_1(u)=\lim_{t_1\to-b}\frac{\phi_0(u)}{1+b^{-1}t_1}\qquad\textrm{and}\qquad
dt'=d\delta_{-b}(t_1)dt_2\cdots dt_n.$$ Note that this limit
exists since there is a factor $1+b^{-1}t_1$ in $c(u^{-1})$ (see
Remark~\ref{rem:A1(iii)}(ii)).

In the exceptional case, let $\pi $ be the measure on
$U=\mathbb{U}\cup\mathbb{U}'$ given by $d\pi
(u)=\frac{W_0(q^{-1})}{|W_0|}\frac{dt}{\phi_0(u)}$ on
$\mathbb{U}$ and $d\pi
(u)=\frac{W_0(q^{-1})}{|W_0'|}\frac{dt'}{\phi_1(u)}$ on
$\mathbb{U}'$, where $W_0'$ is the Coxeter group $C_{n-1}$ (with
$C_1=A_1$ and $C_0=\{1\}$). Then for all $A\in\sca_2$,
\begin{align*}
(A\delta_y)(x)=\int_{U}\widehat{A}(u)\overline{P_{\la}(u)}d\pi(u)\quad\textrm{whenever
$y\in V_{\la}(x)$}
\end{align*}
(see \cite[Theorem~5.7 and Corollary~5.8]{p2}).

To conveniently state formulae in both the standard and exceptional
cases simultaneously, we write $U=\bu$ in the standard case and
(as above) $U=\bu\cup\bu'$ in the exceptional case. Thus (in all
cases), for $A\in\sca_2$,
\begin{align}\label{eq:planinv}
(A\delta_y)(x)=\int_{U}\widehat{A}(u)\overline{P_{\la}(u)}d\pi(u)\quad\textrm{whenever
$y\in V_{\la}(x)$}.
\end{align}

\begin{rem} The form of the Plancherel measure in the exceptional
case requires that $q_1b\geq1$, which follows from a theorem of
D. Higman since the numbers $q_i$, $i\in I$, are the parameters of
a building (see \cite[Lemma~5.6]{p2}). We note that for the
hypergroups associated to the $BC_n$ case the Plancherel measure
is supported on $U=\bu\cup\bu'\cup\bu''\cup\cdots$, where there
are $k$ components, with $k$ defined by
$q_1^{k-1}b\geq1>q_1^{k-2}b$. See \cite[Theorem~5.2.10]{macsph}.
\end{rem}

\section{The Local Limit Theorem}\label{section:2}

The basic approach for the local limit theorem is as follows. Let
$A$ be the transition operator for an isotropic random walk with
matrix $(p(x,y))_{x,y\in V_P}$, as in~(\ref{Adefn}). Then
\begin{align}\label{eq:mk}
p^{(k)}(x,y)=(A^k\delta_y)(x)\qquad\textrm{for all $x,y\in V_P$
and $k\in\bn$}.
\end{align}
Since $\|A\|\leq 1$, we may regard $A$ as in $\sca_2$ and so
$h_u(A)$, $u\in U$, is defined. Writing $\widehat{A}(u)=h_u(A)$
for $u\in U$, we have $\widehat{A}_{\la}(u)=P_{\la}(u)$ and so
\begin{align}\label{Pdefn}
\widehat{A}(u)=\sum_{\la\in P^+}a_{\la}P_{\la}(u).
\end{align}
By (\ref{eq:mk}) and~(\ref{eq:planinv}) we have
\begin{align}\label{prop}
p^{(k)}(x,y)=\int_U\big(\widehat{A}(u)\big)^k\,\overline{P_{\la}(u)}d\pi
(u)\qquad\textrm{whenever $y\in V_{\la}(x)$,}
\end{align}
and we will prove the local limit theorem by determining the
asymptotic behaviour of the integral in (\ref{prop}) as
$k\to\infty$.

\begin{lem}\label{paper3lem1} Let $\la\in P^+$, $\la\neq 0$, $x\in V_P$, and $y\in
V_{\la}(x)$. Then
\begin{enumerate}
\item[(i)] there exists $z\in V_{\la}(x)\cap
V_{\tilde{\alpha}^{\vee}}(y)$, and
\item[(ii)] with $z$ as in $(\mathrm{i})$, there exists $\o\in\O$ such that
$h(y,z;\o)=\tilde{\alpha}^{\vee}$.
\end{enumerate}
\end{lem}

\begin{proof} Note first that if $c$ and $d$ are distinct
$i$-adjacent chambers, $i\in I_P$, with type~$i$ vertices $u$ and
$v$ respectively, then $v\in V_{\at}(u)$ (and $u\in
V_{\tilde{\alpha}^{\vee}}(v)$). To see this, let $\ca$ be any
apartment containing $c$ and $d$, and let $\psi:\ca\to\S$ be a
type rotating isomorphism such that $\psi(u)=0$ and
$\psi(c)=C_0$. Since $\psi(d)$ is $0$-adjacent to $\psi(c)$ we
have $\psi(d)=s_{\tilde{\alpha};1}(C_0)$, and so
$\psi(v)=s_{\tilde{\alpha};1}(0)=\at$. Thus $v\in V_{\at}(u)$.

Part~(i) now follows exactly as in \cite[Lemma~5.1]{cw}; we
include the proof for completeness. Let $\ca$ be an apartment
containing $x$ and $y$, and let $c_0,c_1,\ldots,c_m$ be a gallery
(that is, a sequence of adjacent chambers with $c_{i-1}\neq c_i$
for $1\leq i\leq m$) with $x\in c_0$ and $y\in c_m$ and $m$
minimal. Let $\pi$ be the panel $c_m\backslash\{y\}$, and let
$c'$ be the chamber of $\ca$ with $c'\neq c_m$ and $\pi\subset
c'$ (so $c'=c_{m-1}$ if $m\geq1$). Let $H$ be the wall of $\ca$
determined by~$\pi$, and let $\ca^+$ be the half apartment of
$\ca$ determined by~$H$ containing $c'$. By thickness there exists
a chamber $d\neq c',c_m$ with $\pi\subset d$, and writing $z$ for
the vertex in $d\backslash\pi$ we have $z\in
V_{\tilde{\alpha}^{\vee}}(y)$ by the above discussion. We now
show that $z\in V_{\la}(x)$. By the proof of
\cite[Lemma~9.4]{ronan} there exists an apartment $\cb$
containing $\ca^+\cup d$. Let $\rho_{\ca,c'}$ be the retraction of
$\scx$ onto $\ca$ with center $c'$ (see \cite[\ts{IV.3}]{brown}),
and so the map $\varphi=\rho_{\ca,c'}|_{\cb}:\cb\to\ca$ is a type
preserving isomorphism with $\varphi(x)=x$ and $\varphi(z)=y$
(since $\varphi(d)=c_m$). Since $y\in V_{\la}(x)$ there exists a
type rotating isomorphism $\psi:\ca\to\S$ with $\psi(x)=0$ and
$\psi(y)=\la$ (see \cite[Proposition~5.6(iii)]{p}), and so the map
$\phi=\psi\circ\varphi:\cb\to\S$ is a type rotating isomorphism
with $\phi(x)=0$ and $\phi(z)=\la$. Thus $z\in V_{\la}(x)$.

Part (ii) is a consequence of the following fact.
Let $u,v\in V_P$ with $v\in V_{\la}(u)$. Then there exists
$\o\in\O$ such that $h(u,v;\o)=\la$. To see this, let $\ca$ be
any apartment containing $u$ and $v$, and let $\psi:\ca\to\S$ be
a type rotating isomorphism such that $\psi(u)=0$ and
$\psi(v)=\la$. Let $\o$ be the class of $\psi^{-1}(\cs_0)$. Since
$\psi^{-1}(\cs_0)=\cs^u(\o)$ and
$\psi^{-1}(\la+\cs_0)=\cs^v(\o)$, we have
$\psi^{-1}(\mu)=v_{\mu}^{u}(\o)=v_{\mu-\la}^v(\o)$ for
sufficiently large $\mu\in P^+$, and so $h(u,v;\o)=\la$.
\end{proof}

Recall that
$\bu=\{u\in\Hom(P,\bc^{\times}):|u^{\la}|=1\textrm{ for all
$\la\in P$}\}$. Let
$$\bu_Q=\{u\in\Hom(P,\bc^{\times})\mid u^{\gamma}=1\textrm{ for
all $\gamma\in Q$}\}.$$ Thus $\bu_Q$ is isomorphic to the dual of
the finite abelian group $P/Q$, and so $\bu_Q\cong P/Q$. Hence
$\bu_Q$ is finite, and $\bu_Q\subset\bu$.

\goodbreak\begin{prop}\label{char2} The set $W_0\tilde{\alpha}^{\vee}$ spans $Q$ over $\bz$. Thus if $u\in\Hom(P,\bc^{\times})$ and $u^{w\tilde{\alpha}^{\vee}}=1$ for all $w\in W_0$, then~$u\in\bu_Q$.
\end{prop}

\begin{proof} Let $Q'$ denote the $\bz$-span of $W_0\tilde{\alpha}^{\vee}$. We
show that $R^{\vee}\subset Q'$, from which it follows that
$Q=Q'$, hence the result. Suppose first that $R$ is reduced, and
let $\beta\in R$. By \cite[VI, \ts1, No.3,
Proposition~11]{bourbaki} all roots of a given length are
conjugate under~$W_0$, and so if $|\beta|=|\tilde{\alpha}|$ then
$\beta^{\vee}\in Q'$. Suppose that $|\beta|\neq|\tilde{\alpha}|$.
Since at most 2 root lengths occur in~$R$ (see
\cite[Lemma~10.4.C]{h2}), and since~$R$ is irreducible, there
exists $v,v'\in W_0$ such that $\lan
v\tilde{\alpha},v'\beta\ran\neq0$ (for otherwise
$W_0\tilde{\alpha}\cup W_0\beta$ is a partition of $R$ into
nonempty pairwise orthogonal sets). Thus $\lan
w\tilde{\alpha},\beta\ran\neq0$, where $w=v'^{-1}v$, and so by
\cite[VI, \ts1, No.8, Proposition~25(iv)]{bourbaki} we have $\lan
w\tilde{\alpha}^{\vee},\beta\ran=1$ or $\lan
w\tilde{\alpha}^{\vee},\beta\ran=-1$, depending on if
$w^{-1}\beta\in R^+$ or $w^{-1}\beta\in R^-$. Since
$s_{\beta}(w\tilde{\alpha}^{\vee})=w\tilde{\alpha}^{\vee}-\lan
w\tilde{\alpha}^{\vee},\beta\ran\beta^{\vee}$ we have
$\beta^{\vee}\in Q'$.

Finally, if $R$ is non-reduced, then $Q(R)=Q(R_1)$ and
$W_0(R)=W_0(R_1)$. Since $\tilde{\alpha}$ is also the highest
root of the reduced root system~$R_1$ (with respect to the
natural base), we have $Q'(R)=Q'(R_1)$, and so $Q'=Q$ in all
cases.
\end{proof}

\begin{rem} Since $\tilde{\alpha}$ is a \textit{long root} \cite[Lemma~10.4.D]{h2}
(with the convention that all roots are called long if there is
only one root length), Proposition~\ref{char2} is true whenever
$\tilde{\alpha}$ is replaced with an arbitrary long root~$\alpha$
(for $W_0\alpha=W_0\tilde{\alpha}$). However for general
$\alpha\in R$ the proposition fails, despite the fact that
$W_0\alpha^{\vee}$ spans~$E$ (by \cite[Lemma~10.4.B]{h2}). For
example let $R$ be the standard $B_2$ root system, and take
$\alpha=e_2$. Then the $\bz$-span of $W_0\alpha^{\vee}$ is
$2\bz^2$, whereas~$Q=\bz^2$.
\end{rem}

As usual, if $u,v\in\Hom(P,\bc^{\times})$, define
$uv\in\Hom(P,\bc^{\times})$ by $(uv)^{\la}=u^{\la}v^{\la}$ for all
$\la\in P$.

\goodbreak\begin{lem}\label{eqcond} Let $u\in\mathbb{U}$ and $\la\in P^+$. Then
$|P_{\la}(u)|\leq P_{\la}(1)$, and equality holds for $\la\neq 0$
if and only if $u\in\bu_Q$. Moreover, if $u_0\in\bu_Q$ then
$P_{\la}(u_0u)=u_0^{\la}P_{\la}(u)$ for all
$u\in\Hom(P,\bc^{\times})$.
\end{lem}

\begin{proof}(cf. \cite[Lemma~5.3]{cw}) Let $x,y\in V_P$ be any vertices with $y\in
V_{\la}(x)$. The inequality is clear from the integral formula
(\ref{eq:twoforms}). Suppose equality holds for some $\la\neq 0$.
Write $f(\o)$ for the integrand in (\ref{eq:twoforms}). Then $f$
is a continuous function on $\O$ and $f(\o)\neq 0$ for all
$\o\in\O$. So
$|\int_{\O}f(\o)d\nu_x(\o)|=\int_{\O}|f(\o)|d\nu_x(\o)$ implies
that $f(\o)/|f(\o)|$ is constant, since $\nu_x(O)>0$ for all
non-empty open sets $O\subset\O$. Thus $u^{h(x,y;\o)}$ takes the
constant value $P_{\la}(u)/P_{\la}(1)$ for all $\o\in\O$. Let $z$
be as in Lemma~\ref{paper3lem1}(i). Since the value of the
integral in (\ref{eq:twoforms}) is unchanged if $y$ is replaced
by $z$, it follows that $u^{h(x,y;\o)}=u^{h(x,z;\o)}$ for all
$\o\in\O$. Choosing $\o\in\O$ as in Lemma~\ref{paper3lem1}(ii)
and using the cocycle relations we have
$u^{\tilde{\alpha}^{\vee}}=u^{h(y,z;\o)}=1$. Furthermore, since
the value of the integral in (\ref{eq:twoforms}) is unchanged if
$u$ is replaced by $wu$ for any $w\in W_0$, then
$u^{w\tilde{\alpha}^{\vee}}=1$ for all $w\in W_0$. It follows
from Proposition~\ref{char2} that $u\in\bu_Q$.

Conversely, if $u_0\in\bu_Q$ and $y\in V_{\la}(x)$, then
$u_0^{h(x,y;\o)}=u_0^{\la}$ for all $\o\in\O$, because
$\la-h(x,y;\o)\in Q$ (see \cite[Theorem~3.4(ii)]{p2}). Thus it
follows from (\ref{eq:twoforms}) that
$P_{\la}(u_0u)=u_0^{\la}P_{\la}(u)$ for all
$u\in\Hom(P,\bc^{\times})$. In particular,
$|P_{\la}(u_0)|=P_{\la}(1)$.
\end{proof}

In the following series of estimates we will write $C$ for a
positive constant, whose value may vary from line to line.

For each $\o\in\O$, $x,y\in V_P$ and $1\leq j\leq n$, define
$h_j(x,y;\o)=\lan h(x,y;\o),\alpha_j\ran$.

\goodbreak\begin{lem}\label{bound} Let $x\in V_P$ and $\la\in P^+$. Then $|h(x,y;\o)|\leq |\la|$ and $|h_j(x,y;\o)|\leq C|\la|$
for all $\o\in\O$, all $y\in V_{\la}(x)$, and all $j=1,\ldots,n$.
\end{lem}

\begin{proof} Recall from \cite[Theorem~3.4(ii)]{p2} that $h(x,y;\o)\in\Pi_{\la}$
for all $\o\in \O$ and $y\in V_{\la}(x)$. By~\cite[(2.6.2)]{m} we
have that $\Pi_{\la}\subset\mathrm{conv}(W_0\la)$ (the usual
convex hull in $E$ here), and since $|w\la|=|\la|$ for all $w\in
W_0$, this implies that $|h(x,y;\o)|\leq|\la|$ for all $\o\in\O$
and for all $y\in V_{\la}(x)$. We have $|\lan
h(x,y;\o),\alpha_j\ran|\leq |h(x,y;\o)||\alpha_j|$, proving the
final claim.
\end{proof}

\begin{rem}\label{rem:graph} There is a natural graph with vertex set $V_P$ and
vertices $x,y\in V_P$ joined by an edge if and only if $y\in
V_{\la_i}(x)$ for some $i\in I_0$. In this graph we have
$d(x,y)=\sum_{i=1}^n\lan\la,\alpha_i\ran$ if $y\in V_{\la}(x)$.
Lemma~\ref{bound} shows that $|h(x,y;\o)|$ and $|h_j(x,y;\o)|$
are bounded by $Cd(x,y)$.
\end{rem}

\begin{notation} Let $\t_1,\ldots,\t_n\in\br$ and write
$\t=\t_1\alpha_1+\cdots+\t_n\alpha_n$ (so $\t\in E$). Write
$e^{i\theta}$ for the element of $\Hom(P,\bc^{\times})$ with
$(e^{i\t})^{\la}=e^{i\lan\la,\theta\ran}$ for all $\la\in P^+$.
With this notation (\ref{eq:twoforms}) gives
\begin{align}\label{integral0}
P_{\la}(e^{i\t})=\int_{\O}r^{h(x,y;\o)}e^{i\lan
h(x,y;\o),\t\ran}d\nu_x(\o)\quad\textrm{for all $y\in
V_{\la}(x)$},
\end{align}
and since $P_{\la}(w^{-1}e^{i\t})=P_{\la}(e^{i\t})$ for all $w\in
W_0$, it follows that
\begin{align}\label{integralinv}
P_{\la}(e^{i\t})=\int_{\O}r^{h(x,y;\o)}e^{i\lan
h(x,y;\o),w\t\ran}d\nu_x(\o)\quad\textrm{for all $w\in W_0$,
$y\in V_{\la}(x)$.}
\end{align}
\end{notation}

\goodbreak
\begin{cor}\label{secondapprox} For all $\la\in P^+$, $P_{\la}(e^{i\t})=P_{\la}(1)(1+E_{\la}(\t))$, where
$|E_{\la}(\t)|\leq |\la||\t|$.
\end{cor}

\begin{proof} We have
$$|P_{\la}(e^{i\t})-P_{\la}(1)|\leq\int_{\O}r^{h(x,y;\o)}|e^{i\lan
h(x,y;\o),\t\ran}-1|d\nu_x(\o),$$ and the result follows from
Lemma~\ref{bound} since $|e^{iz}-1|\leq |z|$ for all $z\in\br$.
\end{proof}

Let $\la\in P^+$ and $y\in V_{\la}(x)$. For each $1\leq j,k\leq n$
define
\begin{align}\label{last:eq}
b_{j,k}^{\la}=\frac{1}{2}\int_{\O}h_j(x,y;\o)h_k(x,y;\o)r^{h(x,y;\o)}d\nu_x(\o).
\end{align}
This is independent of the particular pair $x,y\in V_P$ with $y\in
V_{\la}(x)$, for by (\ref{integral0})
$$\frac{\partial^2}{\partial\t_j\partial\t_k}P_{\la}(e^{i\t})\bigg|_{\t=0}=-\int_{\O}h_j(x,y;\o)h_k(x,y;\o)r^{h(x,y;\o)}d\nu_x(\o).$$
(Indeed any expression
$\int_{\O}p(h_1(x,y;\o),\ldots,h_n(x,y;\o))r^{h(x,y;\o)}d\nu_x(\o)$,
where $p$ is a polynomial, is independent of the particular pair
$x,y\in V_P$ with $y\in V_{\la}(x)$).

\goodbreak\begin{lem}\label{thisuu} Let $\la\in P^+$, and $\t_1,\ldots,\t_n\in\br$, and as usual write
$\t=\t_1\alpha_1+\cdots+\t_n\alpha_n$. Then
\begin{align}\label{thisu}
P_{\la}(e^{i\t})=P_{\la}(1)-\sum_{j,k=1}^nb_{j,k}^{\la}\t_j\t_k+R_{\la}(\t)
\end{align}
where $|R_{\la}(\t)|\leq C|\la|^3|\t|^3 P_{\la}(1)$. Furthermore,
$\sum_{j,k=1}^nb_{j,k}^{\la}\t_j\t_k\geq0$, and when $\la\neq 0$,
equality holds if and only if $\t=0$.
\end{lem}

\begin{proof} For $\varphi\in\br$ we have
$e^{i\varphi}=1+i\varphi-\frac{1}{2}\varphi^2+R(\varphi)$ where
$|R(\varphi)|\leq\frac{1}{6}|\varphi|^3$. Applying this to
$\varphi=\lan h(x,y;\o),\t\ran$ and using (\ref{integral0}) we
have
\begin{align*}
P_{\la}(e^{i\theta})&=P_{\la}(1)+i\int_{\O}\lan
h(x,y;\o),\theta\ran r^{h(x,y;\o)}
d\nu_x(\o)\\
&\quad-\frac{1}{2}\int_{\O}\lan
h(x,y;\o),\theta\ran^2r^{h(x,y;\o)}d\nu_x(\o)+R_{\la}(\theta),
\end{align*}
where $|R_{\la}(\theta)|\leq\frac{1}{6}|\lan
h(x,y;\o),\theta\ran|^3P_{\la}(1)\leq\frac{1}{6}|h(x,y;\o)|^3|\theta|^3P_{\la}(1)$.
The bound for $|R_{\la}(\t)|$ follows from Lemma~\ref{bound}.

We claim that for all $j=1,\ldots,n$ and for all $y\in V_P$,
$$\int_{\O}h_j(x,y;\o)r^{h(x,y;\o)}d\nu_x(\o)=0.$$
To see this, let $j\in\{1,\ldots,n\}$ and set $\t=\t_j\alpha_j$
(that is, $\t_k=0$ for all $k\neq j$). By differentiating
(\ref{integralinv}) with respect to $\t_j$, and then evaluating
at $\t_j=0$, firstly with $w=1$ and secondly with $w=s_{j}$, we
see that
$$\int_{\O}h_j(x,y;\o)r^{h(x,y;\o)}d\nu_x(\o)=-\int_{\O}h_j(x,y;\o)r^{h(x,y;\o)}d\nu_x(\o),$$
proving the claim. It is now clear that (\ref{thisu}) holds, and
that $\sum_{j,k=1}^nb_{j,k}^{\la}\t_j\t_k\geq0$. If equality
holds, then
$$\int_{\O}\lan h(x,y;\o),\t\ran^2r^{h(x,y;\o)}d\nu_x(\o)=0.$$
Thus $\lan h(x,y;\o),\t\ran=0$ for almost all $\o\in\O$, and thus
for all $\o\in\O$. Thus, since $\lan h(x,y;\o),t\t\ran=0$ for all
$t\in\br$ and $\o\in\O$, we have $P_{\la}(e^{i(t\t)})=P_{\la}(1)$
for all $t\in\br$ by (\ref{integral0}), and so
$e^{i(t\t)}\in\bu_Q$ for all $t\in\br$ by Lemma~\ref{eqcond}.
Thus $\t=0$ since $|\bu_Q|<\infty$.
\end{proof}

\goodbreak\begin{lem}\label{later:refertothis} There exists a polynomial $p(x_1,\ldots,x_n)$ of
degree at most $M$ such that
\begin{align}\label{eq:y}
P_{\la}(1)=r^{-\la}p(\lan\la,\alpha_1\ran,\ldots,\lan\la,\alpha_n\ran)
\end{align}
for all $\la\in P^+$, where $M>0$ is some integer depending only
on the underlying root system. Furthermore, (by thickness) there exists some
$q>1$ such that
\begin{align}\label{polynomial}
P_{\la}(1)\leq C(|\la|+1)^Mq^{-|\la|}.
\end{align}
\end{lem}

\begin{proof} Assuming that $u^{-\alpha^{\vee}}\neq1$ for any $\alpha\in R_2^+$,
by (\ref{eq:mac}) and the definition of the numbers~$\tau_{\alpha}$ we have
\begin{align}\label{eq:alternatemac2}
c(u)=\prod_{\alpha\in R_2^+}
\frac{(1-\tau_{2\alpha}^{-1}\tau_{\vphantom{2}\alpha}^{-1/2}u^{-\alpha^{\vee}/2})(1+\tau_{\vphantom{2}\alpha}^{-1/2}u^{-\alpha^{\vee}/2})}{1-u^{-\alpha^{\vee}}}.
\end{align}

Write $\s=\la_1+\cdots+\la_n$. It follows from \cite[VI, \ts3,
No.3, Proposition~2]{bourbaki} that
$$\prod_{\alpha\in
R_2^+}(1-u^{-w\alpha^{\vee}})=(-1)^{\ell(w)}u^{\s-w\s}\prod_{\alpha\in
R_2^+}(1-u^{-\alpha^{\vee}})$$ for all $w\in W_0$, and so by
(\ref{eq:mac}) and (\ref{eq:alternatemac2}) we have
\begin{align}\label{eq:lh}
P_{\la}(u)&=r^{-\la}\frac{F(\la)}{\prod_{\alpha\in
R_2^+}(1-u^{-\alpha^{\vee}})}
\end{align}
where $F(\la)$ equals $\frac{1}{W_0(q^{-1})}$ times
\begin{align*}
\sum_{w\in
W_0}\bigg\{(-1)^{\ell(w)}u^{w\la+w\s-\s}\prod_{\alpha\in
R_2^+}(1-\tau_{2\alpha}^{-1}\tau_{\vphantom{2}\alpha}^{-1/2}u^{-w\alpha^{\vee}/2})(1+\tau_{\vphantom{2}\alpha}^{-1/2}u^{-w\alpha^{\vee}/2})\bigg\}.
\end{align*}
We know that $P_{\la}(u)$ is a Laurent polynomial in
$u_1,\ldots,u_n$, and so (\ref{eq:y}) follows from (\ref{eq:lh})
by repeated applications of L'H\^{o}pital's rule. The
inequality~(\ref{polynomial}) follows from
Proposition~\ref{prop:misc}(ii) and the proof of
Proposition~\ref{prop:misc}(iv) in Section~\ref{section:3}.
\end{proof}

Let $A$ be as in (\ref{Adefn}) and $\widehat{A}(u)=h_u(A)$ be as
in (\ref{Pdefn}). It follows from Lemma~\ref{bound} that
$|b_{j,k}^{\la}|\leq C|\la|^2P_{\la}(1)$, and thus the inequality
(\ref{polynomial}) implies that $\sum_{\la\in
P^+}a_{\la}b_{j,k}^{\la}$ is absolutely convergent for each
$1\leq j,k\leq n$. We define
\begin{align}\label{last2:eq}
b_{j,k}=\frac{1}{\widehat{A}(1)}\sum_{\la\in
P^+}a_{\la}b_{j,k}^{\la}.
\end{align}

\begin{cor}\label{uT} Let $A$ be as in (\ref{Adefn}), and let
$\t_1,\ldots,\t_n\in\br$. Then
$$\widehat{A}(e^{i\t})=\widehat{A}(1)\bigg(1-\sum_{j,k=1}^nb_{j,k}\t_j\t_k+R(\t)\bigg),$$
where $\sum_{j,k=1}^nb_{j,k}\t_j\t_k>0$ unless $\t=0$, and where
$|R(\t)|\leq C|\t|^3$.
\end{cor}

\begin{proof} This follows from Lemma~\ref{thisuu}, using (\ref{polynomial}) to bound $R(\t)$.
\end{proof}

\goodbreak\begin{lem}\label{thirdapprox} Let $\t_1,\ldots,\t_n\in\br$. Then
$$\frac{1}{|c(e^{i\t})|^2}=\prod_{\alpha\in
R_2^+}\frac{\lan\alpha^{\vee},\t\ran^2}{\big(1-\tau_{2\alpha}^{-1}\tau_{\vphantom{2}\alpha}^{-1/2}\big)^2\big(1+\tau_{\vphantom{2}\alpha}^{-1/2}\big)^2}(1+E_{\alpha}(\t))$$
where $|E_{\alpha}(\t)|\leq C\lan\alpha^{\vee},\t\ran^2$ for each
$\alpha\in R_2^+$.
\end{lem}

\begin{proof} Observe that for $x\in\br$ and $p>1$
\begin{align}\label{eq:ap1}
\left|\frac{1-e^{-i
x}}{1-p^{-1}e^{-ix}}\right|^2=\frac{x^2}{(1-p^{-1})^2}(1+E_1(x)),
\end{align}
where $|E_1(x)|\leq C x^2$, and for $p>0$
\begin{align}\label{eq:ap2}
\left|\frac{1+e^{-ix}}{1+p^{-1}e^{-ix}}\right|^2=\frac{4}{(1+p^{-1})^2}(1+E_2(x))
\end{align}
where $|E_2(x)|\leq Cx^2$. The result follows by using
(\ref{eq:alternatemac2}), (\ref{eq:ap1}) and~(\ref{eq:ap2}).
\end{proof}

Let
$\mathbb{U}_A=\{u\in\mathbb{U}:|\widehat{A}(u)|=\widehat{A}(1)\}$.
This set will play a role in the local limit theorem, and in the
conditions for irreducibility and aperiodicity of the random
walk. The following lemma gives a description of $\bu_A$ in terms
of the coefficients $a_{\la}$ appearing in~(\ref{Adefn}).

\goodbreak\begin{lem}\label{UA} We have $\bu_A=\{u\in\bu_Q\mid u^{\mu}=u^{\nu}\textrm{ for
all $\mu,\nu\in P^+$ with $a_{\mu},a_{\nu}>0$}\}$. If
$u_0\in\mathbb{U}_A$ then
$\widehat{A}(u_0u)=u_0^{\mu}\widehat{A}(u)$ for all
$u\in\Hom(P,\bc^{\times})$ and all $\mu\in P^+$ such that
$a_{\mu}>0$.
\end{lem}

\begin{proof} For $u\in\bu$ we have
\begin{align}\label{series}
|\widehat{A}(u)|=\bigg|\sum_{\la\in
P^+}a_{\la}P_{\la}(u)\bigg|\leq\sum_{\la\in
P^+}a_{\la}|P_{\la}(u)|\leq\sum_{\la\in
P^+}a_{\la}P_{\la}(1)=\widehat{A}(1).
\end{align}
If $u=u_0\in\bu_A$, then since equality must hold in the second
inequality in (\ref{series}) we have $|P_{\la}(u_0)|=P_{\la}(1)$
whenever $a_{\la}>0$. Since we assume that $a_{\la}>0$ for at
least one nonzero $\la\in P^+$ we have $u_0\in\bu_Q$ by
Lemma~\ref{eqcond}. Thus by Lemma~\ref{eqcond} we have
$P_{\la}(u_0)=u_{0}^{\la}P_{\la}(1)$ for all $\la\in P^+$, and so
since equality must hold in the first inequality in (\ref{series})
we have $u_{0}^{\mu}=u_0^{\nu}$ whenever $a_{\mu},a_{\nu}>0$,
proving that
$$\bu_A\subseteq\{u\in\bu_Q\mid u^{\mu}=u^{\nu}\textrm{ for
all $\mu,\nu\in P^+$ with $a_{\mu},a_{\nu}>0$}\}.$$ Conversely,
if $u_0\in\bu_Q$ and $u_0^{\mu}=u_0^{\nu}$ for all $\mu,\nu\in
P^+$ with $a_{\mu},a_{\nu}>0$, then by Lemma~\ref{eqcond} we see
that $\widehat{A}(u_0u)=u_0^{\mu}\widehat{A}(u)$ for all $u\in\bu$
and any $\mu\in P^+$ such that $a_{\mu}>0$, and so taking $u=1$ we
have $|\widehat{A}(u_0)|=\widehat{A}(1)$, so $u_0\in\bu_A$.
\end{proof}

For $k\in\bn$ and $\la\in P^+$ let
\begin{align*}
I_{k,\la}=\int_{\bu}\big(\widehat{A}(u)\big)^k\,\overline{P_{\la}(u)}d\pi
(u).
\end{align*}
If $y\in V_{\la}(x)$, then by (\ref{prop})
\begin{align}\label{eq:prp}
p^{(k)}(x,y)=\begin{cases}I_{k,\la}&\textrm{in the standard case,
and}\\
I_{k,\la}+I_{k,\la}'&\textrm{in the exceptional case,}\end{cases}
\end{align}
where
\begin{align}\label{eq:prp2}
I_{k,\la}'=\int_{\bu'}\big(\widehat{A}(u)\big)^k\,\overline{P_{\la}(u)}d\pi
(u).
\end{align}
Thus to give an asymptotic formula for $p^{(k)}(x,y)$ we need to
give estimates for $I_{k,\la}$ and~$I_{k,\la}'$.

Given $\epsilon>0$ and $u_0\in\mathbb{U}$, let
$N_{\epsilon}(u_0)=\{u\in\mathbb{U}:|u^{\la_i}-u_0^{\la_i}|<\epsilon\textrm{
for all $i\in I_0$} \}$. Since $|\mathbb{U}_A|<\infty$ we may
choose $\epsilon>0$ sufficiently small so that
\begin{align}\label{eq:epsilon}
N_{\e}(u_0)\cap N_{\e}(u_0')=\emptyset\qquad\textrm{whenever
$u_0,u_0'\in\bu_A$ are distinct}.
\end{align}
Write $N_{\e}=N_{\e}(1)$ and
$N_{\e}(\bu_A)=\bigcup_{u_0\in\bu_A}N_{\e}(u_0)$.

Define $\rho_1=\rho_1(\e)=\sup\{|\widehat{A}(u)|/\widehat{A}(1):
u\in\bu\backslash N_{\e}(\bu_A)\}$, and so $0<\rho_1<1$. Let
$$I_{k,\la}^{\e}=\int_{N_{\e}}\big(\widehat{A}(u)\big)^k\,\overline{P_{\la}(u)}d\pi (u).$$

\goodbreak\begin{lem}\label{thm:1} Fix $\mu\in P^+$ such that $a_{\mu}>0$, and let $\e>0$ satisfy (\ref{eq:epsilon}). If
$u_0^{k\mu}=u_0^{\la}$ for all $u_0\in\bu_A$, then
$$I_{k,\la}=|\bu_A|I_{k,\la}^{\e}+\mathcal{O}\big(\rho_1^{k}\widehat{A}(1)^k\big).$$
Otherwise, $I_{k,\la}=0$.
\end{lem}

\begin{proof} It is clear from the formula for $c(u)$ that $c(u_0u)=c(u)$ for all
$u_0\in\bu_Q$ and $u\in\bu$. Thus by Lemmas~\ref{eqcond}
and~\ref{UA}, if $u_0\in\bu_A$ we have
\begin{align}\label{eq:s}
I_{k,\la}=u_0^{k\mu-\la}\int_{\bu}\big(\widehat{A}(u_0^{-1}u)\big)^k\,\overline{P_{\la}(u_0^{-1}u)}d\pi
(u_0^{-1}u)=u_0^{k\mu-\la}I_{k,\la}.
\end{align}
This shows that $I_{k,\la}=0$ if there exists $u_0\in\bu_A$ such
that $u_0^{k\mu-\la}\neq 1$.

Suppose now that $u_0^{k\mu-\la}=1$ for all $u_0\in\bu_A$. It is
clear that
\begin{align}\label{eq:t}
I_{k,\la}=\int_{N_{\e}(\bu_A)}\big(\widehat{A}(u)\big)^k\,\overline{P_{\la}(u)}d\pi
(u)+\mathcal{O}\big(\rho_1^k\widehat{A}(1)^k\big),
\end{align}
and since $N_{\e}(u_0)=u_0N_{\e}$, the calculation in (\ref{eq:s})
shows that for each $u_0\in\bu_A$,
$$\int_{N_{\e}(u_0)}\big(\widehat{A}(u)\big)^k\,\overline{P_{\la}(u)}d\pi (u)=u_0^{k\mu-\la}\int_{N_{\e}}\big(\widehat{A}(u)\big)^k\overline{P_{\la}(u)}d\pi (u)=I_{k,\la}^{\e},$$
since $u_0^{k\mu-\la}=1$. The result follows from (\ref{eq:t}) by
the choice of~$\e$.
\end{proof}

It is clear from Corollary~\ref{uT} that if each $|\t_j|$,
$j=1,\ldots,n$, is sufficiently small, then
\begin{align}\label{firstapprox}
\widehat{A}(e^{i\t})=\widehat{A}(1)e^{-\sum_{i,j=1}^{n}b_{i,j}\t_i\t_j+G(\t)}\quad\textrm{where}\quad
G(\t)=o\bigg(\sum_{i,j=1}^{n}b_{i,j}\t_i\t_j\bigg).
\end{align}
Writing $\delta=2\sin^{-1}(\e/2)$ we have
$N_{\e}=\{e^{i\t}:|\t_j|<\delta\textrm{ for $j=1,\ldots,n$}\}$,
and so we may choose $\e>0$ sufficiently small so that
\begin{align}\label{firstapproxd}
|G(\t)|\leq\frac{1}{2}\sum_{i,j=1}^nb_{i,j}\t_i\t_j
\end{align}
whenever $e^{i\t}\in N_{\e}$ and $|\t_j|\leq\pi$ for
$j=1,\ldots,n$.

Define constants $K_1$, $K_2$ and $K_3$ by
$K_1=W_0(q^{-1})|W_0|^{-1}(2\pi)^{-n}$,
\begin{align}
\nonumber K_2&=\prod_{\alpha\in
R_2^+}\big(1-\tau_{2\alpha}^{-1}\tau_{\vphantom{2}\alpha}^{-1/2}\big)^{-2}\big(1+\tau_{\vphantom{2}\alpha}^{-1/2}\big)^{-2}\\
\label{eq:integral}K_3&=\int_{\br^n}e^{-\sum_{i,j=1}^{n}b_{i,j}\varphi_i\varphi_j}\prod_{\alpha\in
R_2^+}\lan\alpha^{\vee},\varphi\ran^2d\varphi_1\cdots d\varphi_n,
\end{align}
where $\varphi=\varphi_1\a_1+\cdots+\varphi_n\a_n$.

\goodbreak\begin{lem}\label{thm:2} Let $\e>0$ be such that (\ref{eq:epsilon}) and (\ref{firstapproxd}) hold. Then
$$I_{k,\la}^{\e}=KP_{\la}(1)\widehat{A}(1)^k\,
k^{-|R_2^+|-n/2}\left(1+\mathcal{O}\big(k^{-1/2}\big)\right),$$
where $K=K_1K_2K_3$.
\end{lem}

\begin{proof} Let $\delta=2\sin^{-1}(\e/2)$ as above. We have
\begin{align*}
I_{k,\la}^{\e}&=K_1\int_{-\delta}^{\delta}\cdots\int_{-\delta}^{\delta}\big(\widehat{A}(e^{i\t})\big)^k\,\frac{P_{\la}(e^{-i\t})}{|c(e^{i\t})|^2}\,d\t_1\cdots
d\t_n,
\end{align*}
and so by making the change of variable $\varphi_j=\sqrt{k}\t_j$
for each $j=1,\ldots,n$ we see that
\begin{align}\label{finas}
I_{k,\la}^{\e}&=K_1k^{-n/2}\int_{-\sqrt{k}\delta}^{\sqrt{k}\delta}\cdots\int_{-\sqrt{k}\delta}^{\sqrt{k}\delta}\big(\widehat{A}(e^{i\varphi/\sqrt{k}})\big)^k\,\frac{P_{\la}(e^{-i\varphi/\sqrt{k}})}{|c(e^{i\varphi/\sqrt{k}})|^2}\,d\varphi_1\cdots
d\varphi_n,
\end{align}
where $\varphi=\varphi_1\a_1+\cdots+\varphi_n\a_n$.

By Corollary~\ref{secondapprox} we have
$$P_{\la}(e^{-i\varphi/\sqrt{k}})=P_{\la}(1)(1+E_1(\varphi))\quad\textrm{where}\quad |E_1(\varphi)|\leq\frac{|\la||\varphi|}{\sqrt{k}},$$ and it follows
from Lemma~\ref{thirdapprox} that
$$\frac{1}{|c(e^{i\varphi/\sqrt{k}})|^2}=K_2k^{-|R_2^+|}(1+E_2(\varphi))\prod_{\alpha\in
R_2^+}\lan\alpha^{\vee},\varphi\ran^2,$$ where $|E_2(\varphi)|\leq
k^{-1}p(\varphi_1,\ldots,\varphi_n)$ for some polynomial
$p(x_1,\ldots,x_n)$. Using these estimates (along with
(\ref{firstapprox})) in (\ref{finas}), we see that
$I_{k,\la}^{\e}$ equals
$K_1K_2P_{\la}(1)\widehat{A}(1)^k\,k^{-|R_2^+|-n/2}$ times
\begin{align*}
\int_{X} e^{-\sum_{i,j=1}^n
b_{i,j}\varphi_i\varphi_j+kG(\varphi/\sqrt{k})}\bigg(\prod_{\alpha\in
R_2^+}\lan\alpha^{\vee},\varphi\ran^2\bigg)(1+E_1(\varphi))(1+E_2(\varphi))\,d\varphi_1\cdots
d\varphi_n
\end{align*}
where $X=[-\sqrt{k}\delta,\sqrt{k}\delta]^n$. By (\ref{firstapproxd}), the above integrand is bounded by
$$e^{-\frac{1}{2}\sum_{i,j=1}^nb_{i,j}\varphi_i\varphi_j}\bigg(\prod_{\alpha\in
R_2^+}\lan\alpha^{\vee},\varphi\ran^2\bigg)\bigg(1+\frac{|\la||\varphi|}{\sqrt{k}}\bigg)\bigg(1+\frac{p(\varphi_1,\ldots,\varphi_n)}{k}\bigg),$$
and the lemma follows by the Dominated Convergence Theorem.
\end{proof}

\goodbreak\begin{lem}\label{thm:exrho} Let $\la\in P^+$ and $k\in\bn$. In the exceptional case, there exists
$0<\rho_2<1$ such that
$$
\int_{\bu'}\big(\widehat{A}(u)\big)^k\,\overline{P_{\la}(u)}d\pi
(u)=\mathcal{O}\big(\rho_2^k\widehat{A}(1)^k\big).
$$
\end{lem}

\begin{proof} Let us sketch the proof of this result. The details are given in \cite[Appendix~B.3]{phd}.
Since we are in the exceptional case, we have $R=BC_n$ for some $n\geq1$ and
$q_n<q_0$. Use the isomorphism $\Hom(P,\bc^{\times})\to(\bc^{\times})^n$,
$u\mapsto(t_1,\ldots,t_n)$, where $t_i=u^{e_i}$, to identify
$\Hom(P,\bc^{\times})$ with
$(\bc^{\times})^n$ (and so $\bu$ is identified with $\bt^n$).
Recall that $\bu'$ consists of those $u\in\Hom(P,\bc^{\times})$ such that
$t_1=-\sqrt{q_n/q_0}$ and $t_j\in\bt$ for $2\leq j\leq n$. Write
$\xi_t=(-\sqrt{q_n/q_0},t_2,\ldots,t_n)$, and set $t_j=e^{i\t_j}$ for $2\leq
j\leq n$.

We claim that $|P_{\la}(\xi_t)|<P_{\la}(1)$ for all $\la\neq 0$
and all $t_2,\ldots,t_n\in\bt$, from which the result clearly
follows (since $\widehat{A}$ is continuous on $\bu'$, and $\bu'$
is compact). The first step is to explicitly compute
$P_{\la_1}(u)$ for arbitrary $u\in\Hom(P,\bc^{\times})$.
By~\cite[Lemma~B.3.2]{phd} we have
\begin{align*}
P_{\la_1}(u)=N_{\la_1}^{-1}\left((q_0-1)(1+q_1+\cdots+q_1^{n-1})+\sqrt{q_0q_n}q_1^{n-1}\sum_{j=1}^n(t_j+t_j^{-1})\right)
\end{align*}
for all $u\in\Hom(P,\bc^{\times})$. From this formula we deduce
that $|P_{\la_1}(\xi_t)|<P_{\la_1}(1)$ for all
$t_2,\ldots,t_n\in\bt$ (see \cite[Theorem~B.3.3]{phd}). We now use
this fact to show that $|P_{\la}(\xi_t)|<P_{\la}(1)$ for all
$\la\neq 0$ and all $t_2,\ldots,t_n\in\bt$.

Recall from \cite[Corollary~5.22]{p} that the operators
$\{A_{\la}\}_{\la\in P^+}$ satisfy
$$A_{\la}A_{\mu}=\sum_{\nu\in P^+}a_{\la,\mu;\nu}A_{\nu},$$
where
$$a_{\la,\mu;\nu}=\frac{N_{\nu}}{N_{\la}N_{\mu}}|V_{\la}(x)\cap
V_{\mu^*}(y)|\geq0,$$
and where $x,y\in V_P$ is any pair with $y\in V_{\nu}(x)$. Since
$\la_1=\tilde{\alpha}^{\vee}$ here, an analogous argument to that
given in Lemma~\ref{paper3lem1}(i) shows that $a_{\la,\la;\la_1}>0$ for all
$\la\neq 0$ (see \cite[Lemma~B.3.4]{phd}).

Since the algebra homomorphisms $h_{\xi_t}:\sca\to\bc$ are continuous with
respect to the $\ell^2$-operator norm, and since $\|A_{\mu}\|=P_{\mu}(1)$ for
all $\mu\in P^+$ (see
\cite[Theorem~6.3]{p2}), we have $|P_{\mu}(\xi_t)|\leq P_{\mu}(1)$ for all
$\mu\in P^+$. Hence for $\la\neq0$,
$$|P_{\la}(\xi_t)|^2=|h_{\xi_t}(A_{\la}^2)|\leq\sum_{\mu\in
P^+}a_{\la,\la;\mu}|P_{\mu}(\xi_t)|<\sum_{\mu\in P^+}a_{\la,\la;\mu}P_{\mu}(1)=P_{\la}(1)^2,$$
where we have used the facts that $|P_{\la_1}(\xi_t)|<P_{\la_1}(1)$ and
$a_{\la,\la;\la_1}>0$.
\end{proof}

We now give our local limit theorem.
\begin{thm}\label{cor:llt} Let $y\in V_{\la}(x)$ and $k\in\bn$, and suppose that
$a_{\mu}>0$. If $u_0^{k\mu}=u_0^{\la}$ for all $u_0\in \bu_A$,
then
$$
p^{(k)}(x,y)=|\bu_A|KP_{\la}(1)\widehat{A}(1)^k\,k^{-|R_2^+|-n/2}\big(1+\mathcal{O}(k^{-1/2})\big),$$
where $K$ is as in Lemma~\ref{thm:2}. If $u_0^{k\mu}\neq
u_0^{\la}$ for some $u_0\in\bu_A$, then $p^{(k)}(x,y)=0$.
\end{thm}
\begin{proof} In the standard case the result follows from
(\ref{eq:prp}) and Lemmas~\ref{thm:1} and~\ref{thm:2}. In the
exceptional case, $Q=P$, and so $\bu_Q=\{1\}$, and so
$\bu_A=\{1\}$. The result now follows from (\ref{eq:prp}) and
Lemmas~\ref{thm:1}, \ref{thm:2}, and~\ref{thm:exrho}.
\end{proof}

A random walk on a state-space $X$ is called \textit{irreducible}
if for each pair $x,y\in X$ there exists $k=k(x,y)\in\bn$ such
that $p^{(k)}(x,y)>0$. The \textit{period} of an irreducible
random walk is $\mathfrak{p}=\gcd\{k\geq1\mid p^{(k)}(x,x)>0\}$,
which is independent of $x\in X$ by irreducibility (see
\cite{woess}). An irreducible random walk is called
\textit{aperiodic} if $\mathfrak{p}=1$.

\goodbreak
\begin{cor} Let $A$ be as in (\ref{Adefn}), and suppose that $a_{\mu}>0$. Then
\begin{enumerate}
\item[(i)] $A$ is irreducible if and only if for each $\la\in P^+$
there exists $k=k(\la)\in\bn$ such that $u_0^{k\mu}=u_0^{\la}$
for all $u_0\in\bu_A$, and
\item[(ii)] $A$ is irreducible and aperiodic if and only if
$|\bu_A|=1$.
\end{enumerate}
\end{cor}

\begin{proof} First let us note that in the exeptional case it is easy to see that any walk
with $a_{\mu}>0$ for some $\mu\neq 0$ is both aperiodic and
irreducible, and since $Q=P$ we have $\bu_A=\{1\}$. So consider
the standard case, and suppose that $a_{\mu}>0$. Let $y\in
V_{\la}(x)$. If $A$ is irreducible, then there exists $k\in\bn$
such that $p^{(k)}(x,y)>0$, and so $u_0^{k\mu}=u_0^{\la}$ for all
$u_0\in\bu_A$, by (\ref{eq:prp}) and Lemma~\ref{thm:1}.
Conversely, if for each $\la\in P^+$ there exists $k_0\in\bn$
such that $u_0^{k_0\mu}=u_0^{\la}$ for all $u_0\in\bu_A$, then
writing $r=|\bu_A|$ we have $u_0^{(k_0+rl)\mu}=u_0^{\la}$ for all
$u_0\in\bu_A$ and all $l\geq0$. As $k\to\infty$ through the
values $k_0+rl$, Theorem~\ref{cor:llt} implies irreducibility.

If $|\bu_A|=1$ then $A$ is clearly irreducible, and
Theorem~\ref{cor:llt} shows that $A$ is aperiodic. Conversely, if
$A$ is irreducible and aperiodic, then
\begin{align*}
1=\gcd\{k\geq 1\mid p^{(k)}(x,x)>0\}=\gcd\{k\geq1\mid
u_0^{k\mu}=1\textrm{ for all $u_0\in\bu_A$}\},
\end{align*}
and so $\bu_A=\{1\}$.
\end{proof}

\begin{rem} It is possible to explicitly compute the constant $K_3$ from
(\ref{eq:integral}) (at least in most cases). We refer the reader
to \cite[Section~8.4]{phd} for details. A key step in the
calculation is to observe that there is a number $b>0$ such that
$b_{j,k}=\lan\alpha_j,\alpha_k\ran b$ for all $1\leq j,k\leq n$,
and so $K_3=b^{-|R_2^+|-n/2}J$ where
$$J=\int_{\br^n}e^{-\sum_{j,k=1}^n\lan\alpha_j,\alpha_k\ran\t_j\t_k}\prod_{\alpha\in
R_2^+}\lan\alpha^{\vee},\t\ran^2\,d\t_1\ldots d\t_n$$ and
$\t=\t_1\alpha_1+\cdots+\t_n\alpha_n$. The integral $J$ depends
only on the underlying root system, and has been computed using
Gram's identity in the cases when $R=B_n,C_n,D_n$ or~$BC_n$
(there are other techniques using orthogonal polynomials). We have
\begin{align*}
J=\begin{cases}\pi^{n/2}2^{-n(n-1)}\prod_{i=1}^n(2i)!&\textrm{if $R=B_n$ or $R=BC_n$}\\
\pi^{n/2}2^{-n^2-n-1}\prod_{i=1}^n(2i)!&\textrm{if $R=C_n$}\\
\pi^{n/2}2^{-n^2+n-1}n!\prod_{i=1}^{n-1}(2i)!&\textrm{if $R=D_n$.}
\end{cases}
\end{align*}
When $R=A_n$ the integral $J$ may be written as
$$\int_{\br^n}e^{-(x_1^2+\cdots+x_{n+1}^2)}\prod_{1\leq i<j\leq
n+1}(x_i-x_j)^2\,dx_1\cdots dx_n$$ (up to some constant factors),
where $x_{n+1}=-(x_1+\cdots+x_n)$. We have been unable to compute
this integral. In principle the integrals for the $E,F$ and $G$
cases could be explicitly computed using a computer package.
\end{rem}

\begin{rem}\label{remark:2} Let us briefly discuss some applications of our local limit
theorem to probability measures on groups acting on~$\scx$. An
automorphism $\psi$ of $\scx$ is called \textit{type rotating} if
there exists a type rotating automorphism $\s$ of the Coxeter
graph of~$W$ (in the sense of \cite[\ts4.8]{p}) such that
$\tau(\psi(x))=\s(\tau(x))$ for all $x\in V$. Suppose that $G$ is
a locally compact group acting on $V$ such that $G$ acts
transitively on $V_P$ and such that for each $x\in V_P$ and $g\in
G$ the automorphism $x\mapsto gx$ is type rotating. Assume that
$K=\{g\in G\mid go=o\}$ acts transitively on each set
$V_{\la}(o)$, $\la \in P^+$, where $o\in V_P$ is some fixed
vertex. Let $\varphi$ be the density function of a
bi-$K$-invariant probability measure on~$G$. Then, exactly as in
\cite[Lemma~8.1]{cw}, setting $p(go,ho)=\varphi(g^{-1}h)$ for
$g,h\in G$ defines an isotropic random walk on~$V_P$. Since the
$k$-th convolution power $\varphi^{(*k)}(g)$ is simply
$p^{(k)}(o,go)$, Theorem~\ref{cor:llt} may immediately be
interpreted as a local limit theorem for bi-$K$-invariant
probability measures on~$G$ (the assumption that $a_{\la}>0$ for
some $\la\neq0$ simply means that $\varphi$ is not the indicator
function on~$K$).

As an important modification, suppose now that $G$ is a group of
type preserving simplicial complex automorphisms acting
\textit{strongly transitively} on $\scx$, meaning that~$G$ acts
transitively on pairs $(\ca,c)$ of apartments $\ca$ and chambers
$c\subset\ca$. Fix an apartment $\ca_0$ and a chamber
$c_0\subset\ca_0$. The subgroups $B=\stab_{G}(c_0)$ and
$N=\stab_{G}(\ca_0)$ form a $BN$-pair in $G$ with associated Weyl
group $N/(B\cap N)$ isomorphic to~$W$ \cite[Theorem~5.2]{ronan}.
Indeed the set of left cosets $\{gB\mid g\in G\}$ defines an
affine building (as a chamber system) isomorphic to~$\scx$, where
$gB\sim_i hB$ if and only if $g^{-1}h\in B\lan s_i\ran B$ (where
$wB$ means $nB$ for any $n\in N$ with image $w\in W$). Let $o$ be
the type $0$ vertex of $c_0$. The subgroup $K=\stab_{G}(o)$
of~$G$ equals $BW_0B=\bigcup_{w\in W_0}BwB$ (see
\cite[Theorem~5.4(iii)]{ronan}), and since $G$ acts strongly
transitively and $B\cap N$ is transitive on the chambers of
$\ca_0$, it follows that $K$ is transitive on each set
$V_{\la}(o)$, $\la\in Q\cap P^+$.

Let $\varphi$ be the density function of a bi-$K$-invariant
probability measure on~$G$. To study convolution powers
$\varphi^{(*k)}(g)$, $g\in G$, it is natural to study an
associated random walk on $V_Q=\{x\in V_P\mid\tau(x)=0\}\subseteq
V_P$, where we define $p(go,ho)=\varphi(g^{-1}h)$ for $g,h\in G$.
To apply our local limit theorem we consider these random walks
as reducible isotropic random walks on $V_P$ by setting
$p(x,y)=\varphi(g^{-1}h)$ if $y\in V_{\la}(x)$ and $go\in
V_{\la}(ho)$ for some $\la\in P^+$ (necessarily $\la\in Q\cap
P^+$), and $p(x,y)=0$ otherwise. These random walks have the
property that $p(x,y)=0$ if $\tau(x)\neq \tau(y)$, and it is
simple to see that they are indeed isotropic.
Theorem~\ref{cor:llt} is now applicable, and in particular, by
taking $\scx$ to be the Bruhat-Tits building of a group of
$p$-adic type (see Remark~\ref{rem:A1(ii)}(ii)) we have a local
limit theorem for these groups.

Finally we remark that the methods here can be extended to deal
with groups acting (in a type rotating fashion) on subsets $V_L$
of~$V_P$. Here $L$ is a lattice in~$E$ with $Q\subseteq L\subseteq
P$, and $V_L=\{x\in V_P\mid\tau(x)\in I_L\}$, where
$I_L=\{\tau(\la)\mid\la\in L\}\subseteq I_P$. Thus $V_Q\subseteq
V_L\subseteq V_P$, and our discussion above deals with the extreme
cases of $L=P$ and $L=Q$.
\end{rem}

\section{The Rate of Escape Theorem}\label{section:3}

Let $X$ be any set, and let $P=(p(x,y))_{x,y\in X}$ be a
transition probability matrix. Let $X=\bigcup_{i\in I}X_i$ be a
partition of $X$. We call $P$ \textit{factorisable over $I$} if
for each $i,j\in I$, the sum
$$\sum_{y\in X_j}p(x,y)$$
has the same value for all $x\in X_i$. In this case we write
$\overline{p}(i,j)$ for this value, and let
$\overline{P}=(\overline{p}(i,j))_{i,j\in I}$. Clearly
$\overline{p}(i,j)\geq0$ for all $i,j\in I$, and for each $i\in
I$,
$$\sum_{j\in I}\overline{p}(i,j)=\sum_{j\in I}\sum_{y\in
X_j}p(x,y)=\sum_{y\in X}p(x,y)=1$$ where $x\in X_i$. Thus
$\overline{P}$ is a transition probability matrix (on $I$).
Furthermore, if $(Z_k)_{k\in\bn}$ is a Markov chain on $X$ with
transition probability matrix $P$, then
$(\overline{Z}_k)_{k\in\bn}$, where $\overline{Z}_k=i$ if $Z_k\in
X_i$, defines a Markov chain on $I$ with transition probability
matrix $\overline{P}$.

In our setting, consider the partition (for fixed $o\in V_P$ and
$\o\in \O$) $V_P=\bigcup_{\la\in P}V_{\la}$, where
$$V_{\la}=\{x\in V_P\mid h(o,x;\o)=\la\}.$$

\begin{prop}\label{prop:fact} The matrices (operators) $A_{\la}=(p_{\la}(x,y))_{x,y\in V_P}$, $\la\in P^+$, are
factorisable over~$P$. Moreover, $\overline{p}_{\la}(\mu,\nu)$
does not depend on $o$ or $\o$, and
$\overline{p}_{\la}(\mu,\nu)=\overline{p}_{\la}(0,\nu-\mu)$.
\end{prop}

\begin{proof} Let
$\mu,\nu\in P$ and $x\in V_{\mu}$. By the cocycle relations we
have $h(o,y;\o)=h(x,y;\o)+\mu$ for all $y\in V_P$, and so
\begin{align}\label{eq:starlate}\begin{aligned}
\sum_{y\in V_{\nu}}p_{\la}(x,y)&=\frac{1}{N_{\la}}|\{y\in
V_{\la}(x)\mid h(o,y;\o)=\nu\}|\\
&=\frac{1}{N_{\la}}|\{y\in V_{\la}(x)\mid h(x,y;\o)=\nu-\mu\}|.
\end{aligned}
\end{align}
It follows from \cite[Lemma~3.19]{p2} that $A_{\la}$ is
factorisable, and that $\overline{p}_{\la}(\mu,\nu)$ does not depend
on~$\o\in\O$ or~$o\in V_P$. The transitional invariance is clear.
\end{proof}

\begin{cor}\label{cor:fact} Let $A=(p(x,y))_{x,y\in V_P}$ be as in (\ref{Adefn}).
Then $A$ is factorisable over~$P$. Moreover, for each $\mu,\nu\in
P$ we have $\overline{p}(\mu,\nu)=\overline{p}(0,\nu-\mu)$, and
this value does not depend on $o$ and $\o$. Finally, if
$(Z_k)_{k\in\bn}$ is a Markov chain with transition probability
matrix $A$, then $\overline{Z}_k=h(o,Z_k;\o)$, so that
$\overline{p}(\mu,\nu)=\mathbb{P}(h(o,Z_{k+1};\o)=\nu\mid
h(o,Z_k;\o)=\mu)$.
\end{cor}

\begin{proof} The first statements follow easily from
Proposition~\ref{prop:fact} and the elementary fact that a
(finite or infinite) convex combination of factorisable
transition matrices is again factorisable. The final claim is
immediate from the definition of $\overline{Z}_k$.
\end{proof}

Let $\{T_j\}_{j\in J}$ be a partition of $R_2$ according to root
length (so $|J|=1$ or $2$). For $j\in J$, let $T_j^+=R_2^+\cap
T_j$, and $B_j=B\cap T_j$ (so $B=\bigcup_{j\in J}B_j$, as
$B\subset R_2$). For each $j\in J$, let
$$\rho_j=\frac{1}{2}\sum_{\alpha\in T_j^+}\alpha.$$
Finally, for each $j\in J$ fix some $\beta_j\in T_j^+$.
\goodbreak
\begin{prop}\label{prop:misc} With the above definitions:
\begin{enumerate}
\item[(i)] For $i,j\in J$, if $\alpha\in B_i$, then
$\lan\alpha^{\vee},\rho_j\ran=\delta_{i,j}$. Thus, for each $j\in
J$, $\rho_j\in \overline{\cs}_0$.
\item[(ii)] Let $\la\in P$. Then
$$r^{\la}=\prod_{j\in
J}(\tau_{\beta_j}\tau_{2\beta_j}^2)^{\lan\la,\rho_j\ran}$$
(note that this product has at most two factors).
\item[(iii)] $r^{w_0\la}=r^{-\la}$ for all $\la\in P$ (that is,
$r^{\la^*}=r^{\la}$).
\item[(iv)] If $\la\in P$ and $\mu\preceq\la$, then $r^{\mu}\leq
r^{\la}$, with equality if and only if $\mu=\la$.
\item[(v)] For $w\in W_0$, we have
$c(wr)=\delta_{w,1}W_0(q^{-1})$.
\end{enumerate}
\end{prop}

\begin{proof} (i) Since $\alpha\in B$, $s_{\alpha}$ permutes
$R_2^+\backslash\{\alpha\}$, and since the sets $T_j$, $j\in J$,
are $W_0$-invariant, we have $s_{\alpha}(T_j^+)=T_j^+$ if $j\in
J\backslash\{i\}$. Thus for any $j\in J$ we have
$$s_{\alpha}(\rho_j)=\rho_j-\delta_{i,j}\alpha,$$
and so $\lan\alpha^{\vee},\rho_j\ran=\delta_{i,j}$. Then
$\lan\alpha^{\vee},\rho_j\ran\geq0$ for all $\alpha\in B$, and so
$\rho_j\in\overline{\cs}_0$ for all $j\in J$.

(ii) By (\ref{eq:r}) and the fact that $R_1^+\backslash
R_3^+=2(R_2^+\backslash R_3^+)$ we calculate
\begin{align*}
r^{\la}=\bigg(\prod_{\alpha\in
R_3^+}\tau_{\alpha}^{\frac{1}{2}\lan\la,\alpha\ran}\bigg)\bigg(\prod_{\beta\in
R_2^+\backslash
R_3^+}(\tau_{\beta}\tau_{2\beta}^2)^{\frac{1}{2}\lan\la,\beta\ran}\bigg).
\end{align*}
Since $\tau_{2\alpha}=1$ if $\alpha\in R_3^+$, and since
$\tau_{\beta}=\tau_{\beta_j}$ if $\beta\in T_j$, it follows that
\begin{align*}
r^{\la}&=\prod_{\beta\in
R_2^+}(\tau_{\beta}\tau_{2\beta}^2)^{\frac{1}{2}\lan\la,\beta\ran}=\prod_{j\in
J}(\tau_{\beta_j}\tau_{2\beta_j}^2)^{\lan\la,\rho_j\ran}.
\end{align*}

(iii) Since $w_0\rho_j=-\rho_j$ for $j\in J$, by (ii) we have
$r^{w_0\la}=r^{-\la}$ for all $\la\in P$.

(iv) Observe that $\tau_{\alpha}\tau_{2\alpha}^2=q_{\alpha}$ if
$\alpha\in R_3$, and
$\tau_{\alpha}\tau_{2\alpha}^2=q_0q_{\alpha}$ if $\alpha\in
R_2\backslash R_3$. Thus, by thickness,
$\tau_{\alpha}\tau_{2\alpha}^2>1$ for all $\alpha\in R_2$. Since
$\rho_j\in\overline{\cs}_0$ for $j\in J$, and since
$\mu\preceq\la$ implies that $\la-\mu\in Q^+$, it follows
from~(ii) that $r^{\la-\mu}\geq 1$, with equality if and only if
$\mu=\la$ (for if $\mu\neq\la$, then $\lan\la-\mu,\rho_j\ran>0$
for at least one $j\in J$).

(v) Observe that if $w\neq 1$, then
$wR_2^+\cap(-B)\neq\emptyset$. To see this, if $\alpha\in R_2^+$,
and if $-\alpha\notin wR_2^+$, then $-w^{-1}\alpha\notin R_2^+$,
and so $w^{-1}\alpha\in R_2^+$. It follows that if
$wR_2^+\cap(-B)=\emptyset$, then $w^{-1}B\subset R_2^+$, and so
$w^{-1}R_2^+=R_2^+$. Thus $w=1$ (for by \cite[VI, \ts1, No.6,
Corollary~2]{bourbaki} we have $\ell(w)=|\{\alpha\in R_2^+\mid
w\alpha\in R_2^{-}\}|$ for all $w\in W_0$).

Suppose that $w\neq 1$, and take $\alpha\in R_2^+$ such that
$w\alpha=-\beta\in-B$. Then by~(i) and~(ii),
$$r^{-w\alpha^{\vee}}=r^{\beta^{\vee}}=\tau_{\beta}\tau_{2\beta}^2=\tau_{\alpha}\tau_{2\alpha}^2,$$
and so
$1-\tau_{2\alpha}^{-1}\tau_{\alpha}^{-1/2}r^{-w\alpha^{\vee}/2}=0$.
Thus by (\ref{eq:alternatemac2}) we see
that $c(wr)=0$ whenever~$w\neq1$.
Since $h_u:\sca\to\bc$ is a non-trivial algebra homomorphism we have
\begin{align*}
1&=h_{u}(A_{0})=\frac{1}{W_0(q^{-1})}\sum_{w\in W_0}c(wu)
\end{align*}
for all nonsingular $u\in\Hom(P,\bc^{\times})$. Evaluating at $u=r$ shows that
~$c(r)=W_0(q^{-1})$.
\end{proof}

\begin{rem} Most of Proposition~\ref{prop:misc} can be found on
page~61 of \cite{macsph}. Notice in particular that
Proposition~\ref{prop:misc}(v)
gives a nice factorisation of the Poincar\'{e}
polynomial of~$W_0$, namely $W_0(q^{-1})=c(r)$ (at least when the $q$'s come from a building). See also~\cite{mac3}.
\end{rem}

For each $j\in I_0$, let $j^*\in I_0$ be defined by
$-w_0\alpha_j=\alpha_{j^*}$. Note that $(j^*)^*=j$.

\begin{cor}\label{cor:rform} Let $x\in V_P$, $\la\in P^+$, and $y\in
V_{\la}(x)$.
\begin{enumerate}
\item[(i)] We have
\begin{align*}
P_{\la}(re^{i\t})=\int_{\O}e^{i\lan h(x,y;\o),w_0\t\ran}d\nu_x(\o)
\end{align*}
where $w_0$ is the longest element of $W_0$.
\item[(ii)] For each $j\in I_0$ the integral
$$\gamma_j^{(\la)}=\int_{\O}h_{j^*}(y,x;\o)d\nu_x(\o)$$
is independent of the particular pair $x,y\in V_P$ with $y\in
V_{\la}(x)$ (the $j^*$ here makes the statements of the main
theorems simpler).
\end{enumerate}
\end{cor}
\begin{proof} By
Proposition~\ref{prop:misc}(iii) (and the fact that
$w_0^{-1}=w_0$) we have
\begin{align*}
P_{\la}(re^{i\t})&=P_{\la}(w_0(re^{i\t}))=P_{\la}(r^{-1}w_0(e^{i\t}))=\int_{\O}e^{i\lan
h(x,y;\o),w_0\t\ran}d\nu_x(\o),
\end{align*}
proving (i).

Since $w_0\t=-\sum_{j=1}^n\t_j\alpha_{j^*}$, by~(i) we have
\begin{align}\label{eq:intf}
\frac{\partial}{\partial\t_j}P_{\la}(re^{i\t})\big|_{\t=0}=i\gamma_j^{(\la)},
\end{align}
proving (ii).
\end{proof}

The following proposition gives a symmetry property of the
numbers $\gamma_{j}^{(\la)}$ generalising
\cite[Proposition~3.5(iii)]{cw}. We will not use this result in
this paper.
\begin{prop} Let $j\in I_0$ and $\la\in P^+$.
We have $\gamma_j^{(\la^*)}=\gamma_{j^*}^{(\la)}$.
\end{prop}

\begin{proof} Observe that for $u\in\Hom(P,\bc^{\times})$ and $\la\in P^+$,
$P_{\la^*}(u)=P_{\la}(u^{-1})$. It suffices to prove this for
nonsingular~$u\in\Hom(P,\bc^{\times})$. Using Proposition~\ref{prop:misc}(iii)
we see that if $u\in \Hom(P,\bc^{\times})$ is nonsingular then
$$c(w_0u)=\prod_{\alpha\in
R^+}\frac{1-\tau_{\alpha}^{-1}\tau_{\alpha/2}^{-1/2}u^{-w_0\alpha^{\vee}}}{1-\tau_{\alpha/2}^{-1/2}u^{-w_0\alpha^{\vee}}}=\prod_{\alpha\in
R^+}\frac{1-\tau_{\alpha}^{-1}\tau_{\alpha/2}^{-1/2}u^{\alpha^{\vee}}}{1-\tau_{\alpha/2}^{-\alpha/2}u^{\alpha^{\vee}}}=c(u^{-1})$$
(we have used the facts that $w_0R^+=R^-$ and
$\tau_{\alpha}=\tau_{\beta}$ if $\beta\in W_0\alpha$). Thus
\begin{align*}
P_{\la^*}(u)=\frac{r^{-\la^*}}{W_0(q^{-1})}\sum_{w\in
W_0}c(wu)u^{w\la^*}=\frac{r^{-\la}}{W_0(q^{-1})}\sum_{w\in
W_0}c(ww_0u)u^{ww_0\la^*}=P_{\la}(u^{-1}),
\end{align*}
and so by (\ref{eq:intf}) we have
\begin{align*}
\gamma_j^{(\la^*)}&=-i\frac{\partial}{\partial\t_j}P_{\la}(r^{-1}e^{-i\t})\big|_{\t=0}=-i\frac{\partial}{\partial\t_j}P_{\la}(w_0(re^{-iw_0\t}))\big|_{\t=0}=-i\frac{\partial}{\partial\varphi_{j^*}}P_{\la}(re^{i\varphi})\big|_{\varphi=0},
\end{align*}
and so $\gamma_j^{(\la^*)}=\gamma_{j^*}^{(\la)}$.
\end{proof}

\begin{lem}\label{lem:est} Let $\la\in P^+$ and $j\in I_0$. Then
$$\gamma_j^{(\la)}=\lan\la,\alpha_j\ran+\mathcal{O}(1).$$
\end{lem}
\begin{proof} Let us temporarily write $\rho_{\t}$ in place of $re^{i\t}$. Then
for $w\in W_0$,
$$\frac{\partial}{\partial \t_j}c(w\rho_{\t})\rho_{\t}^{w\la}=i\lan
w\la,\alpha_j\ran
c(w\rho_{\t})\rho_{\t}^{w\la}+\rho_{\t}^{w\la}\frac{\partial}{\partial
\t_j}c(w\rho_{\t}).$$ It follows from Proposition~\ref{prop:misc}(v) that
for $w\in W_0$,
\begin{align*}
\frac{r^{-\la}}{W_0(q^{-1})}\,\frac{\partial}{\partial
\t_j}c(w\rho_{\t})\rho_{\t}^{w\la}\bigg|_{\t=0}=\begin{cases}i\lan
\la,\alpha_j\ran+\mathcal{O}(1)&\textrm{if $w=1$}\\
\mathcal{O}(1)&\textrm{if $w\neq 1$.}\end{cases}
\end{align*}
The result follows from (\ref{eq:intf}).
\end{proof}

\begin{lem} If $(Z_k)_{k\in\bn}$ is a Markov chain in $V_P$ with
$Z_0=x$ and transition operator $A_{\la}$, then for any
$\o\in\O$, $\mathbb{E}(h_j(Z_1,x;\o))=\gamma_{j^*}^{(\la)}$.
\end{lem}

\begin{proof} Since $Z_1\in V_{\la}(x)$ with probability 1, we
have $\int_{\O}h_j(Z_1,x;\o)d\nu_x(\o)=\gamma_{j^*}^{(\la)}$. As
in \cite[Proposition~3.5(ii)]{cw} we see that we may take
expectations under the integral sign, and so
$\gamma_{j^*}^{(\la)}=\int_{\O}\mathbb{E}(h_j(Z_1,x;\o))d\nu_x(\o)$.
By Corollary~\ref{cor:fact}, the distribution of $h_j(Z_1,x;\o)$,
and hence $\mathbb{E}(h_j(Z_1,x;\o))$, is independent of
$\o\in\O$. The result follows.
\end{proof}

We now prove our rate of escape theorem.

\begin{thm}\label{thm:roe} Let $A$ be as in (\ref{Adefn}), and suppose that $\sum_{\la\in
P^+}|\la|a_{\la}<\infty$. Let $(Z_k)_{k\in\bn}$ be the
corresponding Markov chain, and for each $k\in\bn$ let $\nu_k\in
P^+$ be such that $Z_k\in V_{\nu_k}(x)$, where $x=Z_0$. Then for each $j\in I_0$,
with probability 1
$$\frac{1}{k}\lan\nu_k,\alpha_j\ran\to\gamma_j\quad\textrm{as
$k\to\infty$,}$$ where $\gamma_j=\sum_{\la\in
P^+}a_{\la}\gamma_j^{(\la)}$. That is,
$\frac{1}{k}\nu_k\to\gamma_1\la_1+\cdots+\gamma_n\la_n$.
\end{thm}

\begin{proof} Observe first that $\gamma_j<\infty$ by Lemma~\ref{bound} and the
finite first moment assumption. By
Lemma~\ref{lem:est} we have
$\frac{1}{k}\lan\nu_k,\alpha_j\ran=\frac{1}{k}\gamma_j^{(\nu_k)}+\mathcal{O}(k^{-1})$,
and so it suffices to prove that
$$\int_{\O}\frac{h_{j^*}(Z_k,x;\o)}{k}d\nu_x(\o)\to\gamma_j$$
with probability~1.

By Corollary~\ref{cor:fact} we see that for each fixed $\o\in\O$,
$h_{j^*}(Z_k,x;\o)$ is a random variable distributed like a sum
of $k$ independent real random variables, each with the
distribution of $h_{j^*}(Z_1,x;\o)$. Now
$\mathbb{E}(h_{j^*}(Z_1,x;\o))=\gamma_j$, and so by the classical
law of large numbers we have
$$\frac{h_{j^*}(Z_k,x;\o)}{k}\to\gamma_j$$
with probability $1$.

By Remark~\ref{rem:graph} and the second part of
\cite[Proposition~8.8(a)]{woess} we see that
$h_{j^*}(Z_k,x;\o)/k$ is bounded with probability~1. Thus by the
Bounded Convergence Theorem we have
$$\lim_{k\to\infty}\int_{\O}\frac{h_{j^*}(Z_k,x;\o)}{k}d\nu_x(\o)=\int_{\O}\lim_{k\to\infty}\frac{h_{j^*}(Z_k,x;\o)}{k}d\nu_x(\o)=\gamma_j$$
with probability~1, completing the proof.
\end{proof}

\begin{cor}\label{cor:nonneg} The numbers $\gamma_j$, $j=1,\ldots,n$, from Theorem~\ref{thm:roe}
are nonnegative.
\end{cor}

\begin{proof} This is immediate from the rate of escape theorem, since $\nu_k\in
P^+$ for each $k\in\bn$.
\end{proof}

We can strengthen Corollary~\ref{cor:nonneg} by applying the local
limit theorem.

\begin{thm}\label{thm:pos} For $j=1,\ldots,n$ we have $\gamma_j>0$ (where, as always, we assume
that $a_{\la}\neq 0$ for at least one $\la\neq0$). In particular, by taking
$A=A_{\la}$ we have $\gamma_j^{(\la)}>0$ for all $\la\neq 0$. 
\end{thm}

\begin{proof} We follow the outline given in~\cite[Remark~4.7]{cw}. By our
local limit theorem we may choose a pair $(\nu,K)\in P^+\times\bn$ with $K$
large and each $\lan \nu,\alpha_j\ran$, $j=1,\ldots,n$, large, such that
$p^{(K)}(x,y)>0$ whenever $y\in V_{\nu}(x)$. With $A$ as in~(\ref{Adefn}), write
$A^K=\sum_{\la\in P^+}a_{\la}^{(K)}A_{\la}$, and so $a_{\nu}^{(K)}>0$. For each $j=1,\ldots,n$ let
\begin{align}
\label{eq:a1}\gamma_{j,K}&=\sum_{\la\in
P^+}a_{\la}^{(K)}\gamma_j^{(\la)}\\
\label{eq:a2}&=-i\frac{\partial}{\partial\t_j}\widehat{A}^K(re^{i\t})\big|_{\t=0},
\end{align}
where we have used~(\ref{eq:intf}). Note that $\gamma_{j,K}$ is simply the
$\gamma_j$ for the transition matrix $A^K$.

It follows from~(\ref{eq:a2}) and~(\ref{eq:intf}) that $\gamma_{j,K}=K\gamma_j$,
and from Corollary~\ref{cor:nonneg} and~(\ref{eq:a1}) we have $\gamma_{j,K}\geq
a_{\nu}^{(K)}\gamma_j^{(\nu)}$. Now, by Lemma~\ref{lem:est} we see that each
$\gamma_j^{(\nu)}$, $j=1,\ldots,n$, is strictly positive (remember that each
component of $\nu$ may be chosen to be large), and thus $\gamma_{j,K}\geq
a_{\nu}^{(K)}\gamma_j^{(\nu)}>0$. Thus $\gamma_j=\frac{1}{K}\gamma_{j,K}>0$ for
each $j=1,\ldots,n$.
\end{proof}

\begin{rem} Fix a vertex $o\in V_P$, and
recall that $\cs^o(\o)$ denotes the unique sector in the
class~$\o$ based at~$o$, and that for each $\la\in P^+$ we write
$v_{\la}^o(\o)$ for the unique vertex in $\cs^o(\o)\cap
V_{\la}(o)$. Given vertices $x,y\in V_P$, there is a natural
notion of the \textit{convex hull} $\mathrm{conv}\{x,y\}$, as
studied in \cite[Appendix~B]{p2}. We say that a sequence
$(x_k)_{k\in\bn}$ of vertices in $V_P$ \textit{converges to
$\o\in\O$} if for each $\la\in P^+$ there exists $k_{\la}\in\bn$
such that $v_{\la}^o(\o)$ is in $\mathrm{conv}\{o,x_{k}\}$
whenever $k\geq k_{\la}$. It is easy to see that this definition
is independent of the $o\in V_P$ chosen. Theorem~\ref{thm:pos}
shows that for an isotropic random walk $(Z_k)_{k\in\bn}$ we
have, with probability~1, $Z_k\to\o$ for some random
element~$\o\in\O$. The key point to observe to show this is that if $Z_k\in
V_{\nu_k}(Z_0)$, then by Theorems~\ref{thm:roe} and~\ref{thm:pos}
$\lan\nu_k,\alpha_j\ran$, $j=1,\ldots,n$, becomes large as
$k\to\infty$.
\end{rem}

\begin{rem}\label{rem:rw}
We note that the random walk $(\overline{Z}_k)_{k\in\bn}$ on~$P$
from Corollary~\ref{cor:fact} can be explicitly studied using
classical methods, since $P\cong\bz^n$. In the notation
of~(\ref{eq:late2}), by~(\ref{eq:starlate}) and \cite[Lemma~3.19
and~Theorem~6.2]{p2} we have
$\overline{p}_{\la}(0,\mu)=r^{-\mu}a_{\la,\mu}$. Assuming that
$\sum_{\mu\in P}|\mu|\,\overline{p}(0,\mu)<\infty$, a calculation
using~(\ref{eq:late2}) and~(\ref{eq:intf}) shows that the mean
$\mathfrak{m}=\sum_{\mu\in P}\mu \,\overline{p}(0,\mu)$ of the
random walk $(\overline{Z}_k)_{k\in\bn}$ is
$\mathfrak{m}=\sum_{j=1}^n\gamma_{j^*}\la_j$, where $\gamma_j$ is
as in Theorem~\ref{thm:roe}. A similar calculation shows that the
characteristic function for the walk is
$$\sum_{\mu\in
P}\overline{p}(0,\mu)e^{i\lan\mu,\t\ran}=\widehat{A}(r^{-1}e^{i\t}).$$
By Corollary~\ref{cor:fact} this walk is transitionally
invariant, and so the usual Fourier inversion (as in
\cite[\ts{II}.6, Proposition~3]{spitzer}) gives
$$\overline{p}^{(k)}(0,\mu)=\frac{1}{(2\pi)^{n}}\int_{-\pi}^{\pi}\cdots\int_{-\pi}^{\pi}\big(\widehat{A}(r^{-1}e^{i\t})\big)^ke^{-i\lan\mu,\t\ran}\,d\t_1\cdots
d\t_n.$$ The asymptotic behaviour may now be extracted using the
methods in \cite[\ts{III.13}]{woess} and the calculations in
Lemma~\ref{lem:posdef}.
\end{rem}

\section{The Central Limit Theorem}\label{section:4}

\begin{lem}\label{lem:cltlem1} Let $\la\in P^+$. There exists a constant $C$, independent of $\t$ and
$\la$, such that
$$\big|h_{re^{i\t}}(A_{\la})-e^{i\lan\la,\theta\ran}\big|\leq
C|\t|.$$
\end{lem}

\begin{proof} Recall that $r^{w\la}\leq r^{\la}$ and
$c(wr)=\delta_{w,1}W_0(q^{-1})$ for all
$w\in W_0$ (see Proposition~\ref{prop:misc}). Thus
\begin{align*}
\big|h_{re^{i\t}}(A_{\la})-e^{i\lan\la,\t\ran}\big|\leq\frac{1}{W_0(q^{-1})}\sum_{w\in
W_0}\big|c(w(re^{i\t}))-c(wr)\big|.
\end{align*}
The result follows since each $c(w(re^{i\t}))$ is a smooth
function in $\t_1,\ldots,\t_n$.
\end{proof}

\begin{lem}(See Theorem~\ref{boundedop})
\label{lem:cltlem3} The homomorphisms
$h_{re^{i\t}}:\sca\to\bc$, $\t\in E$, are bounded.
\end{lem}

\begin{proof} For each $\la\in P^+$ we have
$|h_{re^{i\t}}(A_{\la})|\leq 1$ by Corollary~\ref{cor:rform}(i).
\end{proof}

Let $x\in V_P$. The \textit{spherical function (with respect to
$x$)} associated to $h_u$ is the function $F_u^x:V_P\to\bc$ which
for each $\la\in P^+$ takes the constant value $h_u(A_{\la})$ on
the set $V_{\la}(x)$.

\begin{lem}\label{lem:cltlem4} Let $Z_0=x$, and suppose that $u\in\Hom(P,\bc^{\times})$ is such that
$h_u:\sca\to\bc$ is bounded. Then
$\mathbb{E}(F_u^x(Z_k))=(\widehat{A}(u))^k$.
\end{lem}

\begin{proof} We have $A^k=\sum_{\la\in
P^+}a_{\la}^{(k)}A_{\la}$ where $a_{\la}^{(k)}=\mathbb{P}(Z_k\in
V_{\la}(x))$. Since $F_u^x(Z_k)=h_u(A_{\la})$ if $Z_k\in
V_{\la}(x)$, we have
$$\mathbb{E}(F_u^x(Z_k))=\sum_{\la\in
P^+}a_{\la}^{(k)}h_u(A_{\la})=h_u(A^k)=(\widehat{A}(u))^k,$$
where we have used the continuity of $h_u$ on the closure of
$\sca$ in the space of bounded linear operators on $\ell^1(V_P)$
to justify the last two equalities.
\end{proof}

For $1\leq j,k\leq n$ and $\la\in P^+$, let
\begin{align}\label{gjk}
\gamma_{j,k}^{(\la)}=\int_{\O}h_{j^*}(y,x;\o)h_{k^*}(y,x;\o)d\nu_x(\o)
\end{align}
where $x,y\in V_P$ are any vertices with $y\in V_{\la}(x)$ (as in
Corollary~\ref{cor:rform}(ii) this is easily seen to be
independent of the particular pair $x,y\in V_P$ with $y\in
V_{\la}(x)$ chosen).
If we suppose a finite second moment assumption:
\begin{align}\label{fsm}
\sum_{\la\in P^+}|\la|^2a_{\la}<\infty,
\end{align}
then for all $1\leq j,k\leq n$ we have $\sum_{\la\in
P^+}a_{\la}\gamma_{j,k}^{(\la)}<\infty$, and we denote this value
by $\gamma_{j,k}$.

\begin{lem}\label{lem:posdef} Suppose that (\ref{fsm}) holds. Then with
$\gamma_j$, $1\leq j\leq n$ as defined in Theorem~\ref{thm:roe},
and $\gamma_{j,k}$, $1\leq j,k\leq n$ as defined above,
$$\widehat{A}(re^{i\t})=1+i\sum_{j=1}^n\gamma_j\t_j-\frac{1}{2}\sum_{j,k=1}^n\gamma_{j,k}\t_j\t_k+o(|\t|^2).$$
Furthermore, if $\t\neq 0$ then
$\left(\sum_{j=1}^n\gamma_j\t_j\right)^2<\sum_{j,k=1}^n\gamma_{j,k}\t_j\t_k$.
\end{lem}

\begin{proof} Consider the case $A=A_{\la}$, $\la\neq 0$. Using Corollary~\ref{cor:rform}(i), the
elementary result
$e^{i\varphi}=1+i\varphi-\frac{1}{2}\varphi^2+o(\varphi^2)$
implies that
$$\widehat{A}_{\la}(re^{i\t})=1+i\sum_{j=1}^n\gamma_j^{(\la)}\t_j-\frac{1}{2}\sum_{j,k=1}^n\gamma_{j,k}^{(\la)}\t_j\t_k+o(|\la||\t|),$$
where we have used Lemma~\ref{bound}. The first claim follows.

To deduce the final claim, let
$B_{\la}=\sum_{j=1}^n\gamma_j^{(\la)}\t_j$ and
$C_{\la}=\sum_{j,k=1}^n\gamma_{j,k}^{(\la)}\t_j\t_k$. Then
$$B_{\la}^2=\bigg(\int_{\O}\sum_{j=1}^{n}h_{j^*}(y,x;\o)d\nu_x(\o)\bigg)^2\leq\int_{\O}\bigg(\sum_{j=1}^n
h_{j^*}(y,x;\o)\bigg)^2d\nu_x(\o)=C_{\la},$$ and
$$\bigg(\sum_{j=1}^n\gamma_j\t_j\bigg)^2=\bigg(\sum_{\la\in
P^+}a_{\la}B_{\la}\bigg)^2\leq\sum_{\la\in P^+}
a_{\la}B_{\la}^2\leq\sum_{\la\in
P^+}a_{\la}C_{\la}=\sum_{j,k=1}^n\gamma_{j,k}\t_j\t_k.$$

To see that the inequality is strict if $\t\neq 0$, recall that by
hypothesis there exists $\la\neq 0$ such that $a_{\la}>0$. If
equality holds in the inequality $B_{\la}^2\leq C_{\la}$, then
for $y\in V_{\la}(x)$,
$$\lan h(x,y;\o),w_0\t\ran=\sum_{j=1}^n h_{j^*}(y,x;\o)\t_j$$
is independent of $\o\in \O$, and thus by
Corollary~\ref{cor:rform}(ii) this quantity is independent of the
particular pair $x,y\in V_P$ with $y\in V_{\la}(x)$ too. Choosing
$z\in V_{\la}(x)\cap V_{\tilde{\alpha}^{\vee}}(y)$ as in
Lemma~\ref{paper3lem1}(i), we have
$$\lan h(y,z;\o),w_0\t\ran=\lan h(x,z;\o),w_0\t\ran-\lan h(x,y;\o),w_0\t\ran=0.$$
By modifying the proof of Lemma~\ref{paper3lem1}(ii), it is easy
to see that for each $w\in W_0$ there exists $\o_w\in\O$ such
that $h(y,z;\o_w)=w\tilde{\alpha}^{\vee}$, and thus by the above
$\lan w\tilde{\alpha}^{\vee},w_0\t\ran=0$ for all $w\in W_0$.
Thus $\t=0$, since $W_0\tilde{\alpha}^{\vee}$ spans~$E$
\cite[Lemma~10.4.B]{h2}.
\end{proof}

Let $\G_1(\t)=\sum_{j=1}^n\gamma_j\t_j$ and
$\G_2(\t)=\sum_{j,k=1}^n\gamma_{j,k}\t_j\t_k$, and write
$$\G(\t)=\G_2(\t)-\G_1^2(\t)=\sum_{j,k=1}^ng_{j,k}\t_j\t_k.$$
By Lemma~\ref{lem:posdef}, $\G=(g_{j,k})_{j,k=1}^n$ is a positive
definite matrix.

\begin{thm} Let $A$ be as in (\ref{Adefn}) and suppose that
(\ref{fsm}) holds. As in Theorem~\ref{thm:roe}, for each
$k\in\bn$ let $\nu_k\in P^+$ be such that $Z_k\in V_{\nu_k}(x)$,
where $x=Z_0$. Then
$$\left(\frac{\lan\nu_k,\alpha_1\ran-\gamma_1k}{\sqrt{k}},\ldots,\frac{\lan\nu_k,\alpha_n\ran-\gamma_nk}{\sqrt{k}}\right)$$
converges in distribution to the normal distribution $N(0,\G)$,
with $\G$ as above.
\end{thm}

\begin{proof} Following the proof of the classical Central Limit Theorem (see \cite[Proposition~II.8]{spitzer} for example), it suffices to
show that
\begin{align}\label{eq:cltequiv}
\lim_{k\to\infty}\mathbb{E}(e^{i(\lan\nu_k,\t\ran-k\G_1(\t))/\sqrt{k}})=e^{-\frac{1}{2}\G(\t)}.
\end{align}

By Lemma~\ref{lem:cltlem1} we have
$$e^{i\lan\nu_k,\t\ran/\sqrt{k}}=P_{\nu_k}(re^{i\t/\sqrt{k}})+o(k^{-1/2})=F_{re^{i\t/\sqrt{k}}}^x(Z_k)+o(k^{-1/2}),$$
and so by Lemmas~\ref{lem:cltlem3} and~\ref{lem:cltlem4} we have
$$\mathbb{E}(e^{i\lan\nu_k,\t\ran/\sqrt{k}})=(\widehat{A}(re^{i\t/\sqrt{k}}))^k+o(k^{-1/2}).$$
Thus using Lemma~\ref{lem:posdef} we compute
\begin{align*}
\mathbb{E}(e^{i(\lan\nu_k,\t\ran-k\G_1(\t))/\sqrt{k}})&=\left(1-\frac{1}{2k}\G(\t)+o(k^{-1})\right)^k+o(k^{-1/2})=e^{-\frac{1}{2}\G(\t)}+o(k^{-1/2}).
\end{align*}
Thus (\ref{eq:cltequiv}) holds, completing the proof.
\end{proof}

\begin{appendix}

\section{Bounded Spherical Functions}\label{appendix:1}
It is easy to see that each $A\in\sca$ maps $\ell^1(V_P)$ into
itself. Let $\sca_1$ denote the closure of $\sca$ in the space
$\mathscr{L}(\ell^1(V_P))$ of bounded linear operators on
$\ell^1(V_P)$. Thus $\sca_1$ is a commutative unital Banach
$*$-algebra. The algebra homomorphisms $h:\sca_1\to\bc$ are precisely the
extensions of those algebra homomorphisms $h_u:\sca\to\bc$ which are bounded. In
this appendix we determine the $u\in\Hom(P,\bc^{\times})$ for which this holds.

In the notation of Remark~\ref{rem:A1(ii)}(ii), it is shown in
\cite[Theorem~4.7.1]{macsph} that $h_u:\sca_Q\to\bc$ is bounded if and only
$|u^{w\la}|\leq r^{\la}$ for all $\la\in Q\cap P^+$ and all $w\in W_0$. The
proof given in \cite{macsph} requires some knowledge of the \textit{singular
cases} (when the denominator of a $c(wu)$ function vanishes). While it should
be possible to generalise the proof in \cite{macsph} to cover the more general
setting of homomorphisms $h_u:\sca\to\bc$, we will
provide a different proof which does not require any
specific details of the singular cases (instead our proof uses the Plancherel
measure).

We restrict our attention to the standard case (where
$\tau_{\alpha}\geq1$ for all $\alpha\in R$). In the exceptional
case (where $\tau_{\alpha}<1$ for some $\alpha\in R$) we have
$R=BC_n$ for some $n\geq1$ and so $Q=P$ and $\sca_Q=\sca$. Thus
Macdonald's analysis in \cite{macsph} covers this specific case.

\begin{rem} (i) In fact in \cite[Theorem~4.7.1]{macsph} Macdonald proves a geometric analog of the result stated
above. For $u\in\Hom(Q,\bc^{\times})$, identify
$\log|u|\in\Hom(Q,\br)$ with the unique element $x_u\in E$ which
satisfies $\lan\la,x_u\ran=\log|u^{\la}|$ for all $\la\in Q$. Let
$D=\{x_{wr}\mid w\in W_0\}$. Then \cite[Theorem~4.7.1]{macsph}
says that $h_u:\sca_Q\to\bc$ is bounded if and only if
$x_{u}\in\conv(D)$.

(ii) Note that we have already seen in Lemma~\ref{lem:cltlem3} that the homomorphisms
$h_{re^{i\t}}:\sca\to\bc$ are bounded, and that this was enough information to
prove our central limit theorem. It is, of course, still desirable to have the
much more accurate Theorem~\ref{boundedop} below.
\end{rem}

For $\la\in P^+$ and $u\in\Hom(P,\bc^{\times})$, define the
\textit{monomial symmetric function} $m_{\la}(u)$ by
$$m_{\la}(u)=\sum_{\mu\in W_0\la}u^{\mu},$$
where $W_0\la=\{w\la\mid w\in W_0\}$. By \cite[(6.1)]{p2} there are numbers $a_{\la,\mu}$ such that
\begin{align}\label{eq:late2}
P_{\la}(u)=\sum_{\mu\in
P}a_{\la,\mu}u^{\mu}=\sum_{\mu\preceq\la}a_{\la,\mu}m_{\mu}(u)
\end{align}
(where the second sum is over those $\mu\in P^+$ with $\la-\mu\in
Q^+$), and it follows (using \cite[Theorem~6.11]{p} for example)
that for $\la,\mu\in P^+$ there are numbers $b_{\la,\mu}$ such
that
\begin{align}\label{eq:po}
m_{\la}(u)=\sum_{\mu\preceq\la}b_{\la,\mu}P_{\mu}(u).
\end{align}

\begin{lem}\label{lem:app} Let $\la,\mu\in P^+$ and $\mu\preceq\la$. In the standard case there exists a constant $K>0$ independent of $\la$
and $\mu$ such that $|b_{\la,\mu}|\leq Kr^{\mu}$. Thus there is a constant $C>0$ independent of $\la\in
P^+$ such that
$$\sum_{\mu\preceq\la}|b_{\la,\mu}|\leq C|\Pi_{\la}|r^{\la}.$$
\end{lem}

\begin{proof} Since we assume that we are in the standard case, by \cite[Lemma~6.1]{p2} we have
$$b_{\la,\mu}=N_{\mu}\int_{\mathbb{U}} m_{\la}(u)\overline{P_{\mu}(u)}d\pi (u).$$

Using (\ref{eq:N}) and the techniques used to derive \cite[(5.2)]{p2} we see
that
\begin{align*}
b_{\la,\mu}&=\frac{W_0(q^{-1})}{W_{0\mu}(q^{-1})|W_0|}r^{\mu}\int_{\mathbb{U}}\sum_{w\in
W_0}\frac{m_{\la}(wu)u^{-w\mu}\overline{c(wu)}}{|c(wu)|^2}du\\
&=\frac{W_0(q^{-1})}{W_{0\mu}(q^{-1})}r^{\mu}\int_{\mathbb{U}}m_{\la}(u)\frac{u^{-\mu}}{c(u)}du,
\end{align*}
and the first claim easily follows. The final claim follows from Proposition~\ref{prop:misc}(iv).
\end{proof}

Let
$\Upsilon=\{u\in\Hom(P,\bc^{\times}):|u^{w\la}|\leq r^{\la}\textrm{ for
all $\la\in P^+$ and all $w\in W_0$}\}$.

\begin{thm}\label{boundedop} The algebra homomorphism $h_u:\sca\to\bc$ is bounded if and only if $u\in \Upsilon$.
\end{thm}

\begin{proof} In the exceptional case this follows from
\cite[Theorem~4.7.1]{macsph}, as remarked at the beginning of
this appendix. Suppose we are in the standard case. If
$u\in\Upsilon$ is nonsingular, then by (\ref{eq:mac}) we
have
$$|h_{u}(A_{\la})|\leq\frac{1}{W_0(q^{-1})}\sum_{w\in W_0}|c(wu)|\quad\textrm{for all $\la\in P^+$}.$$
Thus $h_{u}:\sca\to\bc$ is bounded, and so by
\cite[Theorem~I.2.5]{davidson}, $|h_{u}(A_{\la})|\leq 1$ for all
$\la\in P^+$.

If $u\in\Upsilon$ is singular, it is clear that there
exists a sequence $(u_{(k)})_{k\in\bn}$ in $\Upsilon$ such that
each $u_{(k)}$ is nonsingular and $u_{(k)}\to u$. By the above we
have $|h_{u_{(k)}}(A_{\la})|\leq 1$ for all $\la\in P^+$, and
since each $h_{u_{(k)}}$ is a Laurent polynomial (in the variables
$\{u_{(k)}^{\la_i}\}_{i\in I_0}$) it follows that
$|h_{u}(A_{\la})|\leq 1$ for all $\la\in P^+$. Thus $h_{u}$ is
bounded for all $u\in \Upsilon$.

Suppose now that $h_{u}:\sca\to\bc$ is bounded (so
$|h_u(A_{\la})|\leq 1$ for all $\la\in P^+$). Then for all
$\la\in P^+$, by (\ref{eq:po}) and Lemma~\ref{lem:app},
$$|m_{\la}(u)|\leq\sum_{\mu\preceq\la}|b_{\la,\mu}|\leq C|\Pi_{\la}|r^{\la}.$$
Thus fixing $\la$ and considering $m_{k\la}(u)$ for $k\in \bn$
gives
$$\bigg|\sum_{\mu\in W_0\la}(r^{-\la}u^{\mu})^k\bigg|\leq p(k),$$
where $p(k)$ is a polynomial. It is elementary that this implies
that $|r^{-\la}u^{\mu}|\leq 1$ for all $\mu\in W_0\la$, hence the
result.
\end{proof}

\end{appendix}

\bibliography{PhD.bib}
\bibliographystyle{plain}

\end{document}